\documentclass[12pt,a4paper,reqno]{amsart}

\usepackage{amssymb}
\usepackage{amsmath}
\usepackage{latexsym}
\usepackage{exscale}
\usepackage[latin1]{inputenc}
\usepackage{epsfig}
\usepackage{verbatim}
\usepackage{graphicx}
\usepackage{subfigure}
\usepackage{epic,eepic}

\newcommand{\norm}[1]{\left\Vert#1\right\Vert}
\newcommand{\abs}[1]{\left| #1 \right|}
\newcommand{\ip}[2]{\langle #1 , #2 \rangle}

\numberwithin{equation}{section} 

\newtheorem{thm}{Theorem}[section] 
\newtheorem{cor}[thm]{Corollary}

\newtheorem{defi}[thm]{Definition}
\newtheorem{rem}[thm]{Remark}

\newtheorem{lem}[thm]{Lemma}

\begin{document}

\allowdisplaybreaks

\title[SHEARLETS, RESTRICTED NONLINEAR APPROXIMATION AND INTERPOLATION]
 {\large DEMOCRACY OF
    SHEARLET FRAMES WITH
    APPLICATIONS TO 
    APPROXIMATION AND INTERPOLATION
    }

\author{Daniel Vera}

\address{Daniel Vera
\\
Departamento de Matem\'aticas
\\
Instituto Tecnológico Aut\'onomo de México
\\
Río Hondo, México. Tel: +52 (55) 5628 4000 ext. 3826, Fax: +52 (55) 5628 4086.}

\email{daniel.vera@itam.mx}


\date{\today}


\begin{abstract} Shearlets on the cone provide Parseval frames
for $L^2$. They also provide near-optimal approximation for the
class $\mathcal{E}$ of cartoon-like images. Moreover, there are spaces
associated to them other than $L^2$ and there exist embeddings
between these and classical spaces.

We prove approximation properties of the cone-adapted shearlets in a
more general context, namely, when the target function belongs to a
class or space different to $\mathcal{E}$ and when the error is not
necessarily measured in the $L^2$-norm but in a much wider family of
smoothness space of high anisotropy. We first prove democracy of
shearlet frames in shear anisotropic inhomogeneous Besov and
Triebel-Lizorkin sequence spaces. Then, we prove embeddings between
approximation spaces and discrete weighted Lorentz spaces in the
framework of shearlets. Simultaneously, we also prove that these
embeddings are equivalent to Jackson and Bernstein type
inequalities. This allows us to find real interpolation between
these highly anisotropic spaces. Finally, we describe how some of
these results are extended to other shearlet and curvelet generated
spaces.
\end{abstract}

\maketitle

\vskip0.5cm
\section{Introduction.}\label{S:Intro}
Wavelets are orthonormal bases for $L^2$ and, more generally, frames
for several function and distribution spaces that have been
successfully applied in harmonic analysis, numerical and analytical
solution of certain partial differential equations, signal
processing and statistical estimation. Not only can they be used to
characterize some classical spaces but there are also implementable
fast algorithms. Although wavelets provide better approximation
properties than Fourier techniques, they lack of high directional
sensitivity in dimensions $d\geq 2$ since the number of wavelets
(from a multi resolution analysis) remain constant across scales:
$2^d-1$. Some directional systems with increasing anisotropy
have been created to overcome this
limitation. Two of them are the \emph{curvelets} of Candès and
Donoho (\cite{CaDo00} and \cite{CaDo04}) and the \emph{shearlets} on
the cone of Guo, Kutyniok and Labate (\cite{GKL06}).

\subsection{Approximation and wavelets}\label{sS:Intro_Appx-Wvlts}Approximation theory benefits from the unconditional bases
provided by wavelet theory since ``it is enough to threshold the
properly normalized wavelet coefficients" (see \cite{DeV}) to
achieve good $N$-term nonlinear approximation. It is also well known
that the approximation order is closely related to the smoothness of
the function. A generalization of the nonlinear approximation
theory, called restricted nonlinear approximation (RNLA), was
carried out by Cohen, DeVore and Hochmuth in \cite{CDH} in the
setting of wavelet bases in Hardy and Besov spaces where the authors
control the measure of the index set of the approximation elements
instead of the number of terms in the approximation. This measure
(generally other than the counting measure) is closely related to a
weighting of the coefficients in the wavelet expansion. One of the
novelties in \cite{CDH} is that the approximation spaces are not
necessarily contained in the space in which the error is measured.
Some extensions of the restricted nonlinear approximation were done
in \cite{KP06} and \cite{HV11} for general bases in quasi-Banach
distribution spaces and quasi-Banach lattice sequence spaces,
respectively. We develop our results based in \cite{HV11} since its
framework is more general than \cite{KP06} (in fact, \cite{HV11}
allows to recover results in \cite{CDH}, in contrast to
\cite{KP06}). Approximation theory is related to real interpolation
 in a way that some identifications
between approximation spaces, interpolations spaces, discrete
Lorentz spaces and Besov spaces can be proved for certain
parameters. These identifications turn out to be equivalent to what
is known as democracy and to Jackson and Bernstein type inequalities
(see \cite{CDH}, \cite{KP06}, \cite{HV11} and references therein).
The concept of democracy (of a basis or a frame in a space) is
related to that of $p$-space which, in turn, is better known as
$p$-Temlyakov property: If, for $C>0$, the condition
\begin{equation}\label{e:p-Temlkv}
\frac{1}{C}(\nu(\Gamma))^{1/p}
    \leq \norm{\sum_{I\in\Gamma}\frac{\mathbf{e}_I}{u_I}}_{\mathfrak{f}}
    \leq C(\nu(\Gamma))^{1/p},
\end{equation}
holds for all $\Gamma\subset\mathcal{D}$ such that
$\nu(\Gamma)<\infty$ we say that $(\mathfrak{f}, \mathcal{E}, \nu)$
(where $\mathfrak{f}$ is a quasi-Banach space, $\mathcal{E}$ is an
unconditional basis and $\nu$ is a measure on the index set)
shares the $p$-Temlyakov property.

The importance of democracy is twofold. First, a basis is greedy
(the error of the -impractical- best $N$-term approximation and the
error of the greedy algorithm are comparable) if and only if it is
unconditional and democratic. This result is due to Konyagin and
Temlyakov \cite{KoTe99}. Second, it has been established in
\cite{CDH}, \cite{GH04}, \cite{KP06} and \cite{HV11}, among others,
several implications regarding equivalence of democracy, embeddings
between approximation spaces and Lorentz spaces and the Jackson and
Bernstein type inequalities. Admissible wavelet bases are
unconditional bases for many classical spaces as Triebel-Lizorkin
and Besov spaces. However, examples of non-democratic bases are
admissible wavelet bases in Besov spaces when $p \neq q$ and Orlicz
spaces when they do not coincide with Lebesgue spaces. Examples of
democratic bases are admissible wavelet bases in classic
Tribel-Lizorkin spaces, among other spaces (with weights).

We now discuss the importance of the Jackson and Bernstein
inequalities and their relation with Approximation Spaces and
Interpolation Spaces. This paragraph is entirely based on Section 3
of \cite{CDH}. Let $\mathbb{X},\mathbb{Y}$ be a pair of spaces
embedded in a Hausdorff space $\mathfrak{X}$. Let
$\mathbb{X}+\mathbb{Y}$ be the space which consists of all functions
$f$ that can be written as $f=g+h$ for $h\in\mathbb{X}$ and
$g\in\mathbb{Y}$. Define the norm on $\mathbb{X}+\mathbb{Y}$ by
$$\norm{f}:= \inf_{f=h+g}\{\norm{h}_{\mathbb{X}} + \norm{g}_{\mathbb{Y}}\}.$$
More generally, Peetre's $K$-functional is defined, for any $t>0$, by
$$K(f,t):=K(f,t,\mathbb{X},\mathbb{Y}):=\inf_{f=g+h}\{\norm{h}_{\mathbb{X}}+t\norm{g}_{\mathbb{Y}}\}.$$
Let the approximation error be defined by
$$\sigma(f,t)_{\mathbb{X}}:=\inf_{g\in\Sigma_t}\norm{f-g}_{\mathbb{X}}.$$
``The usual setting for approximation takes $t=n, n=1,2,\ldots$ and $\mathbb{Y}\subset\mathbb{X}$ but the results are the same (and the proofs almost identical) in this more general setting (Section 3 of \cite{CDH})" of restricted non-linear approximation.
Assume, additionally, that $\mathbb{X}_t$ is a (possibly non-linear) subspace of $\mathbb{X}+\mathbb{Y}$ and $\mathbb{X}_t\subset\mathbb{X}_u$ if $t\leq u$. The Jackson inequality is defined, for some $r>0$, by
$$\sigma(f,t)_{\mathbb{X}}\leq Ct^{-r}\norm{f}_{\mathbb{Y}}, \;\;f\in\mathbb{Y}, \;\;t>0.$$
In order to compare $\sigma$ and $K$, let $\epsilon>0$ arbitrary and $f=f-g+g$ be such that
$$\norm{f-g}_{\mathbb{X}}+ t^{-r}\norm{g}_{\mathbb{Y}} = K(f,t^{-r})+\epsilon.$$
If $S$ is a best approximation to $g$ from $\Sigma_t$ then, from the Jackson inequality,
\begin{eqnarray*}
\nonumber
\sigma(f,t)_{\mathbb{X}}
    &\leq& \norm{f-S}_{\mathbb{X}} \leq \norm{f-g}_{\mathbb{X}} + \norm{g-S}_{\mathbb{X}}\\
    &\leq& K(f,t^{-r})+\epsilon+C t^{-r}\norm{g}_{\mathbb{Y}}\leq C K(f,t^{-r})+\epsilon.
\end{eqnarray*}
Since $\epsilon$ is arbitrary, then
$$\sigma(f,t)_{\mathbb{X}}\leq C K(f, t^{-r}).$$
The Bernstein inequality provides a weak inverse to the previous
inequality. Since the approximation spaces are defined by $\sigma$
and the interpolation spaces by $K$ both spaces can be characterized
by each other when the the Jackson and Bernstein inequalities hold.

\subsection{Shearlets and first approximation results}\label{sS:Intro_Appx-Shrlts}The number of
directions in the shearlets on the cone or the curvelet systems
doubles at each (other) scale yielding near-optimal approximation to
the class $\mathcal{E}$ of so-called \emph{cartoon-like images} made
up of $C^2$-functions in $[0,1]^d$ except in $C^2$-discontinuities.
It has been proved that, for $f\in\mathcal{E}(\mathbb{R}^2)$,
\begin{equation}\label{e:intro_AppxErr-Cartoon}
\norm{f-f_N^D}_2^2\approx N^{-2}(\log N)^3, \text{ as }
N\rightarrow\infty,
\end{equation}
where $f_N^D$ stands for the curvelet or shearlet approximation with
$N$ terms. This is an optimal approximation except for the
logarithmic factor. In contrast, wavelet approximation gives an
error decay of only $O(N^{-1})$. When $d=3$, the rate of error
approximation is $O(N^{-1}(\log N)^2)$ for both directional systems,
whilst for wavelets it decreases only as $O(N^{-1/2})$. For the
previous statements see \cite{Don01} and \cite{GuLa11}.

Both directional systems form tight frames. However, the shearlet
systems that come from group operation are based on wavelets with
composite dilations of Guo, Labate, Lim, Weiss and Wilson in
\cite{GLLWW} which take full advantage of the theory of affine
systems. In contrast, the curvelet systems are based on polar
coordinates. The use of the shearing operator allows a natural
transition from the continuous to the discrete setting and an easier
implementable framework since the discrete shearing operator leaves
the integer lattice unchanged. We will give a more thorough
description about shearlets and the spaces we will focus on in
Section \ref{sS:Shlts}.

\subsection{Smoothness spaces of curvelets and
shearlets}\label{sS:Intro_Shrlts-SmthnssSpcs}A natural question is
whether these directional systems can generate/characterize other
function/distribution spaces as in the case of wavelet systems and,
if so, what the relation with classical function/distribution spaces
is. The first question has been answered affirmatively several times
and, moreover, some embeddings between these new spaces and
classical ones have been found, answering the second question. The
theory of decomposition spaces was applied by Borup and Nielsen in
\cite{BoNi07} to develop what can be called the \emph{curvelet
decomposition spaces} or \emph{curvelet smoothness spaces}. Short
after Dahlke, Kutyniok, Steidl and Teschke defined in \cite{DKST}
the \emph{shearlet coorbit spaces} through the theory of coorbit
spaces. More recently, Labate, Mantovani and Negi applied again the
general theory of decomposition spaces to introduce the
\emph{shearlet smoothness spaces} in \cite{LMN2012}. A more
classical approach is done in \cite{Ver13}, where the author follows
the ideas of Frazier and Jawerth in \cite{FJ85} to develop the
\emph{shear anisotropic inhomogeneous Besov spaces} of functions and
sequences. In \cite{LMN2012} the authors proved that the shearlet
smoothness spaces and the curvelet smoothness spaces are equivalent.
This result considers only band-limited generators. In a more recent
paper Grohs and Kutyniok defined in \cite{GrKu12} the
\emph{parabolic molecules} that encompass all the curvelet-like and
shearlet-like generators, \emph{i.e.} not necessarily band-limited and compact support.
It is also shown in \cite{GrKu12} that all
curvelet-like and shearlet-like (sequence) spaces (of Besov type)
are equivalent. This is achieved by showing the almost orthogonality
between any two systems of parabolic molecules. Another consequence
of this almost orthogonality is the meta-theorem that states that
all frame systems based on parabolic scaling posses the same
approximation properties considering sparse-promoting $\ell^p$
spaces, i.e. $0<p\leq 1$. In this paper we consider more general
Lorentz spaces and define explicitly the approximation spaces. All
of the above spaces are related to classical Besov spaces. This can
easily be verified from their definition in $\ell^q(\ell^p)$
(quasi-)norms in sequence spaces ($\ell^q(L^p)$ quasi-norms in
\cite{Ver13}) and from their embedding results. As another family of
spaces generated/characterized by the shearlet system (this time
with $L^p(\ell^q)$ quasi-norms) is the \emph{shear anisotropic
inhomogeneous Triebel-Lizorkin spaces} (of functions and sequences)
developed in \cite{Ver12}, following this time the ideas of Frazier
and Jawerth in \cite{FJ90}. We observe that, in order to develop a
Triebel-Lizorkin type spaces with shearlets via decomposition
spaces, it is necessary to use the tools developed in \cite{Ver12},
namely the Fefferman-Stein-Peetre maximal function with shear
anisotropic dilations and related inequality (see Lemma 4.2.3 in
\cite{Ver12}). We point out next some differences between these
spaces. The shearlet coorbit spaces theory (with continuous
parameters and non-uniform directional information) allows the study
of homogeneous spaces with Banach frames (after discretizing the
representation). The rest of spaces aforementioned are naturally
inhomogeneous by construction but generate quasi-Banach spaces with
uniform directional information. The papers \cite{Ver12} and
\cite{Ver13} are the only ones that also consider function spaces.
In applications the amount of information is limited by the sampling
operation and memory/transmission restrictions and so the
inhomogeneous setting is well-adapted to computational procedures.

\subsection{Contribution and expected impact}\label{sS:Intro_Contribution}
Most of the research on approximation properties with curvelet-like
and shearlet-like generators has been focused on showing
near-optimality of sparse approximation to the class $\mathcal{E}$
of cartoon-like images in $2D$ and $3D$ with the error measured in
the $L^2$-norm. However, the $L^2$-norm does not give necessarily
the best visually faithful approximation. For this and other
discussions (\emph{v.gr.} statistical estimation) on the reasons to
measure the error on norms different than $L^2$ see Section 10 of
\cite{DeV} and references therein.

Here, we study approximation properties of shearlet systems when the
target function belongs to a certain (shear anisotropic) smoothness
space, more general than the class $\mathcal{E}$, and the error is
measured on a different (shear anisotropic) smoothness space. For
example, the class $\mathcal{E}(\mathbb{R}^2)$ is a small subset of
the class of functions of bounded variation $BV(\mathbb{R}^2)$,
which in turn lies between $B^{1,1}_1(\mathbb{R}^2)$ and
$B^{1-\varepsilon,1}_1(\mathbb{R}^2)$ (see \cite{CDDD} for a
discussion on the advantages of working with
$B^{1,1}_1(\mathbb{R}^2)$ instead of $BV(\mathbb{R}^2)$). 
Working with sequence spaces will not be a limitation since, as in
the case of wavelets, we have that shear anisotropic inhomogeneous
Besov and Triebel-Lizorkin distribution spaces are a retract of
their sequence spaces counterparts. Therefore, we can transfer the
results established here to the function spaces setting of
\cite{Ver13} and \cite{Ver12}.

The results contained in this paper relate explicitly the
approximation error decay with the space to which the target
function belongs to, in the same spirit of
(\ref{e:intro_AppxErr-Cartoon}). It also paves the way to the use of
thresholding algorithms to efficiently compress or reduce noise of
an image ``through" a given (smoothness) space other than $L^2$ with
the use of the shearlet systems.

\subsection{Outline}\label{sS:Intro_Outline} In Section
\ref{S:Defs_Notation} we introduce definitions, previous results and
notation of i) the theory of RNLASS and ii) the shearlet system and
related spaces we will be working with. In Section
\ref{S:Upr-Lwr_Tmlkv-prop_T-L_Shrlts} we prove democracy of shearlet
frames in some shear anisotropic inhomogeneous sequence spaces.
Applications to approximation theory and interpolation are presented
in Section \ref{S:RNLApprx_SeqSpcsShrlts}. Finally, a straight
extension of some of our results to other spaces generated by
shearlets or curvelets is given in Section
\ref{S:Extn_Parblc_Molecls} thanks to the results on parabolic
molecules in \cite{GrKu12}. Section \ref{S:Extn_Parblc_Molecls}
contains also a brief discussion regarding the space
$BV(\mathbb{R}^2)$ and shearlets.

\vskip0.5cm
\section{Definitions and notation}\label{S:Defs_Notation}

\subsection{Approximation and Interpolation in sequence
spaces}\label{sS:RNLASS} Part from this section is based on
\cite{HV11}. Denote by $S$ the space of all sequences
$\mathbf{s}=\{s_I\}_{I\in\mathcal{D}}$ of complex numbers indexed by
a countable set $\mathcal{D}$. Denote by
$\mathcal{E}=\{\mathbf{e}_I\}_{I\in\mathcal{D}}$ the canonical basis
of $S$, this is, $\mathbf{e}_I$ is the element of $S$ with entry $1$
at $I$ and $0$ otherwise. Thus, for $\Gamma\subset \mathcal{D}$, any
element of $S$ can be written as $\sum_{I\in\Gamma}
s_I\mathbf{e}_I$.

\begin{defi}\label{d:quasi-Bnch-seq-lattice}
A linear space of sequences $\mathfrak{f}\subset S$ is a \textbf{quasi-Banach
(sequence) lattice} if there is a quasi-norm $\norm{\cdot}_\mathfrak{f}$ in
$\mathfrak{f}$ with respect to which $\mathfrak{f}$ is complete and satisfies:

(a) Monotonicity: if $\mathbf{t}\in \mathfrak{f}$ and
$\abs{s_I}\leq\abs{t_I}$ for all $I\in\mathcal{D}$, then
$\mathbf{s}\in \mathfrak{f}$ and $\norm{\{s_I\}}_\mathfrak{f}\leq \norm{\{t_I\}}_\mathfrak{f}$.

(b) If $\mathbf{s}\in \mathfrak{f}$, then $lim_{n\to \infty} \norm{ s_{I_n}
\mathbf {e}_{I_n}}_\mathfrak{f} =0$, for some enumeration
$\mathfrak{I}=\{I_1, I_2,\dots\}$.
\end{defi}
We will say that a quasi-Banach (sequence) lattice $\mathfrak{f}$ is
\textbf{embedded} in $S$, and write $\mathfrak{f}\hookrightarrow S$ if
\begin{equation}\label{e:def_f-embdd-S}
\lim_{n\rightarrow\infty} \norm{\mathbf{s}^n-\mathbf{s}}_\mathfrak{f}=0
    \Rightarrow \lim_{n\rightarrow\infty} s_I^{(n)}=s_I \;\;\forall
    I\in\mathcal{D}.
\end{equation}
\begin{rem}\label{r:monotncty_implies_unconditional}
When $\mathcal{E}=\{\mathbf{e}_I\}_{I\in\mathcal{D}}$ is a Schauder
basis for $\mathfrak{f}$, condition (a) in Definition
\ref{d:quasi-Bnch-seq-lattice} implies that $\mathcal{E}$ is an
unconditional basis for $\mathfrak{f}$ with constant $C=1$.
\end{rem}

In this paper $\nu$ will denote a positive measure on the discrete
set $\mathcal{D}$ such that $\nu(I)> 0$ for all $I\in\mathcal{D}$.
In the classical $N$-term approximation $\nu$ is the counting
measure (i.e. $\nu(I)=1$ for all $I\in\mathcal{D}$), but more
general measures are used in the restricted non-linear approximation
case. For $\Gamma\subset\mathcal{D}$, $\nu(\Gamma)=\sum_{I\in\Gamma} \nu(I)$.
The measure $\nu$ will be used to control the type (instead of the number) of terms
in the approximation.

\begin{defi}\label{d:Standard_scheme_Rstrct_NonLin_Apprx}
We say that $(\mathfrak{f},\nu)$ is a \textbf{standard scheme} (for
restricted non linear approximation) if

i) $\mathfrak{f}$ is a quasi-Banach (sequence) lattice embedded in $S$.

ii) $\nu$ is a measure on $\mathcal{D}$ as explained in the
   previous paragraph.
\end{defi}

Let $(\mathfrak{f},\nu)$ be a standard scheme. For $t>0$, define
$$\Sigma_{t,\nu}:=\{\mathbf{t}=\sum_{I\in\Gamma}t_I{\mathbf e}_I: \nu(\Gamma)\leq t\}.$$
Notice that $\Sigma_{t,\nu}$ is not linear, but
$\Sigma_{t,\nu}+\Sigma_{t,\nu}\subset\Sigma_{2t,\nu}$.
Given $\mathbf{s}\in S$, the $\mathfrak{f}$-\textbf{error of
approximation} to $\mathbf{s}$ (or $\mathfrak{f}$-risk) by elements
of $\Sigma_{t,\nu}$ is given by
$$\sigma_\nu(t,\mathbf{s})=\sigma_\nu(t,\mathbf{s})_\mathfrak{f}:=\inf_{\mathbf{t}\in\Sigma_{t,\nu}}
    \norm{\mathbf{s}-\mathbf{t}}_\mathfrak{f}.$$
Notice that elements $\mathbf{s}\in S$ not in $\mathfrak{f}$ could have
finite $\mathfrak{f}$-risk since elements of $\Sigma_{t,\nu}$ could have
infinite number of entries.

\begin{defi}\label{d:Rstrct_Apprx_Spcs}
(\textbf{Restricted Approximation Spaces}) Let $(\mathfrak{f},\nu)$ be a
standard scheme.

i) For $0<\xi<\infty$ and $0<\mu<\infty$,
$\mathcal{A}^\xi_\mu(\mathfrak{f},\nu)$ is defined as the set of all
$\mathbf{s}\in S$ such that
\begin{equation}\label{e:def_RstrctApprxSpcs}
\norm{\mathbf{s}}_{\mathcal{A}^\xi_\mu(\mathfrak{f},\nu)}
    :=
    \left(\int_0^\infty[t^\xi\sigma_\nu(t,\mathbf{s})]^\mu\frac{dt}{t}\right)^{1/\mu}<\infty.
\end{equation}

ii) For $0<\xi<\infty$ and $\mu=\infty$,
$\mathcal{A}^\xi_\infty(\mathfrak{f},\nu)$ is defined as the set of all
$\mathbf{s}\in S$ such that
\begin{equation}\label{e:def_RstrctApprxSpcs_infty}
\norm{\mathbf{s}}_{\mathcal{A}^\xi_\infty(\mathfrak{f},\nu)}
    := \sup_{t>0} t^\xi\sigma_\nu(t,\mathbf{s}) <\infty.
\end{equation}
\end{defi}

Notice that the spaces $\mathcal{A}^\xi_\mu(\mathfrak{f},\nu)$ depend on the
canonical basis $\mathcal E$ of $S$. When $\mathfrak{f}$ is understood, we will
write $\mathcal{A}^\xi_\mu(\nu)$ instead of
$\mathcal{A}^\xi_\mu(\mathfrak{f},\nu)$.

\begin{rem}\label{r:replace_def_RstrctApprxSpcs}
If $\mathbf{s}\in \mathfrak{f}$, using $\sigma_\nu(t,\mathbf{s})\leq
\norm{\mathbf{s}}_\mathfrak{f}$, it is easy to see that
(\ref{e:def_RstrctApprxSpcs}) can be replaced by
$\norm{\mathbf{s}}_\mathfrak{f}$ plus the same integral from $1$ to
$\infty$. We need to consider the whole range $0<t<\infty$ since we
do not assume $\mathbf{s}\in \mathfrak{f}$. Similar remark holds for
$\mu=\infty$ in (\ref{e:def_RstrctApprxSpcs_infty}). Nevertheless,
the properties of the restricted non-linear approximation spaces are
the same as the $N$-term approximation spaces (see Section 5 in
\cite{CDH}. \emph{Grosso modo}, instead of counting the dyadic cubes
that contain $x$, one measures the ``volume" of those dyadic cubes.
Thus, there may be infinite cubes containing $x$ whose measure is
finite).
\end{rem}

The \textbf{discrete Lorentz spaces} $\ell^{p,\mu}(\nu)$,
$0<p<\infty$, $0 < \mu \leq \infty$, are defined as the set of all
sequences $\mathbf{s}=\{s_I\}_{I\in\mathcal{D}}\in S$ such that
\begin{eqnarray*}
  \norm{\mathbf{s}}_{\ell^{p,\mu}(\nu)}
    &:=& \left(\int_0^\infty [t^{1/p}\mathbf{s}^\ast_\nu(t)]^\mu
        \frac{dt}{t}\right)^{1/\mu}, 
\end{eqnarray*}
(usual modifications if $\mu = \infty$), where
$\mathbf{s}^\ast_\nu(t)=
    \inf\{\lambda>0: \nu(\{I\in\mathcal{D}:\abs{s_I}>\lambda\})\leq t\}$
is the non-increasing rearrangement of $\mathbf{s}$ with respect to
the measure $\nu$. When $\mu=p$ we have
$\ell^{p,p}(\nu)=\ell^p(\nu)$ and thus
$$\norm{\mathbf{s}}_{\ell^p(\nu)}=\left(\sum_{I\in\mathcal{D}}\abs{s_I}^p\nu(I)\right)^{1/p}.$$
We use a weight sequence $\mathbf{u}=\{u_I\}_{I\in\mathcal{D}}$,
$u_I>0$, to control the weight of each $s_I$ as follows. The spaces
$\ell^{p,\mu}(\mathbf{u},\nu)$ are given by $
\norm{\mathbf{s}}_{\ell^{p,\mu}(\mathbf{u},\nu)}
    :=
    \norm{\{u_Is_I\}_{I\in\mathcal{D}}}_{\ell^{p,\mu}(\nu)}<\infty$.

For $0<\theta<1$ and $0<q\leq \infty$ the \textbf{interpolation space} $(\mathbb{X},\mathbb{Y})_{\theta,q}$ is defined as the set of all functions $f\in \mathbb{X}+\mathbb{Y}$ such that, for $0<q<\infty$,
$$|f|_{(\mathbb{X},\mathbb{Y})_{\theta,q}}:= \left(\int_0^\infty [t^{-\theta}
    K(f,t)]^q \frac{dt}{t}\right)^{1/q}< \infty,$$
with the usual modification when $q=\infty$.

Finally, we denote $A\hookrightarrow B$ the continuous inclusion of
$A$ into $B$.

We will use the next particular cases of Theorems 2.6.1 and 2.7.2 in \cite{HV11}.
\begin{thm}\label{t:J-B_ineq_equiv_for_seq}
Let $(\mathfrak{f},\nu)$ be a standard scheme and let
$\mathbf{u}=\{u_I\}_{I\in\mathcal{D}}$ be a weight sequence. Fix
$\xi>0$ and $\mu\in(0,\infty]$. Then, for any $0<p<\infty$ and $r$
such that $\frac{1}{r}=\frac{1}{p}+\xi$, the following are
equivalent:
\begin{enumerate}
  \item[J1)] The upper democracy holds: There exists $C>0$ such that for all $\Gamma\subset\mathcal{D}$
        with $\nu(\Gamma)<\infty$,
        $$\norm{\sum_{I\in\Gamma}\frac{\mathbf{e}_I}{u_I}}_\mathfrak{f} \leq C (\nu(\Gamma))^{1/p}.$$

  \item[J2)] $\ell^{r,\mu}(\mathbf{u},\nu)\hookrightarrow
        \mathcal{A}^\xi_\mu(\mathfrak{f},\nu)$.

  \item[J3)] The space $\ell^{r,\mu}(\mathbf{u},\nu)$ satisfies the
        Jackson inequality of order $\xi$: There exists $C>0$
        such that
        $$\sigma_\nu(t,\mathbf{s})_\mathfrak{f}\leq Ct^{-\xi}\norm{\mathbf{s}}_{\ell^{r,\mu}(\mathbf{u},\nu)}
            , \;\; \text{ for all } \;\;\mathbf{s}\in\ell^{r,\mu}(\mathbf{u},\nu).$$
\end{enumerate}
With the same setting for $(\mathfrak{f},\tilde{\nu})$, the
following are equivalent:
\begin{enumerate}
  \item[B1)] The lower democracy holds: There exists $C>0$ such that for all $\Gamma\subset\mathcal{D}$
        with $\tilde{\nu}(\Gamma)<\infty$,
        $$\frac{1}{C}(\tilde{\nu}(\Gamma))^{1/p}\leq \norm{\sum_{I\in\Gamma} \frac{\mathbf{e}_I}{u_I}}_\mathfrak{f}.$$

  \item[B2)] $\mathcal{A}^\xi_\mu(\mathfrak{f},\tilde{\nu})\hookrightarrow
        \ell^{r,\mu}(\mathbf{u},\tilde{\nu})$.

  \item[B3)] The space $\ell^{r,\mu}(\mathbf{u},\tilde{\nu})$ satisfies the Bernstein
        inequality of order $\xi$: There exists $C>0$ such that
        $$\norm{\mathbf{s}}_{\ell^{r,\mu}(\mathbf{u},\tilde{\nu})}\leq C t^\xi\norm{\mathbf{s}}_\mathfrak{f}
            \;\; \text{ for all }\;\; \mathbf{s}\in\Sigma_{t,\tilde{\nu}}\cap \mathfrak{f}.$$
\end{enumerate}
\end{thm}

Obviously, when $\nu=\tilde{\nu}$ we will say that the upper and
lower democracy ($p$-Temlyakov property) are the same and have
characterizations and identifications.

It is well known that $N$-term approximation and real interpolation
are interconnected. If the Jackson and Bernstein's inequalities hold
for $\nu=$ counting measure, $N$-term approximation spaces are
characterized in terms of interpolation spaces (see \emph{e.g.}
Theorem 3.1 in \cite{DeP88} or Section 9, Chapter 7 in
\cite{DeL93}). As stated in \cite{CDH} (see also Remark
\ref{r:replace_def_RstrctApprxSpcs}), the same scheme is true for
the RNLA. Below we state two more results we need in this paper. The
proofs are straight-forward modifications (see Section 5 in
\cite{CDH}) of those given in the references cited in the previous
paragraph.

\begin{thm}\label{t:J&Bineq_seq->ApprxSpcs_r_InterpolSpcs}
Let $(\mathfrak{f},\nu)$ be a standard scheme. Suppose that the
quasi-Banach lattice $\mathfrak{g}\subset S$ satisfies the Jackson
and Bernstein's inequalities for some $r>0$. Then, for $0<\xi<r$ and
$0<\mu\leq\infty$ we have
$$\mathcal{A}^\xi_\mu(\mathfrak{f},\nu)=\left(\mathfrak{f},\mathfrak{g}\right)_{\xi/r,\mu}.$$
\end{thm}

It is not difficult to show that the spaces
$\mathcal{A}^r_q(\mathfrak{f},\nu), 0<r<\infty, 0<q\leq\infty$,
satisfy the Jackson and Bernstein's inequalities of order $r$, so
that by Theorem \ref{t:J&Bineq_seq->ApprxSpcs_r_InterpolSpcs},
$$\mathcal{A}^\xi_\mu(\mathfrak{f},\nu)=\left(\mathfrak{f},\mathcal{A}^r_q(\mathfrak{f},\nu)\right)_{\xi/r,\mu},$$
for $0<\xi<r$ and $0<\mu\leq\infty$. From here, and using the
reiteration theorem for real interpolation (\emph{e.g.} Ch. 6, Sec.
7 of \cite{DeL93}) we obtain the following result.

\begin{cor}\label{c:Reit_Thm_4_Rstrct_ApprxSpcs}
Let $0<\alpha_0,\alpha_1<\infty, 0<q,q_0,q_1\leq\infty$ and
$0<\theta<1$. Then,
$$\left(\mathcal{A}^{\alpha_0}_{q_0}(\mathfrak{f},\nu),
\mathcal{A}^{\alpha_1}_{q_1}(\mathfrak{f},\nu)\right)_{\theta,q}
    =\mathcal{A}^\alpha_q(\mathfrak{f},\nu), \;\;\; \alpha=(1-\theta)\alpha_0+\theta\alpha_1$$
for a standard scheme $(\mathfrak{f},\nu)$.
\end{cor}


\vskip0.5cm
\subsection{Shearlets}\label{sS:Shlts}
Let $\mathfrak{d}$ be a coordinate in the plane of frequencies
$\hat{\mathbb{R}}^d$. Let $j\geq 0$ and $k\in\mathbb{Z}^d$ be the scale index
and position, respectively, and $\ell$ be the shear parameter such
that $\ell=(\ell_1,\ldots,\ell_{d-1})$ with $-2^j\leq\ell_i\leq2^j$,
$i=1,\ldots,d-1$.

Define the truncated cone
\begin{equation}\label{e:Cone_Domain}
\mathcal{C}^{(1)}=\{(\xi_1,\ldots,\xi_d)\in\hat{\mathbb{R}}^d:
    \abs{\xi_1}\geq\frac{1}{8}, \abs{\frac{\xi_\mathfrak{d}}{\xi_1}}\leq 1, \mathfrak{d}=2,\ldots,d \}.
\end{equation}
Let $\hat{\psi}_1,\hat{\psi}_2\in C^\infty(\mathbb{R})$ with
$\text{supp }\hat{\psi}_1\subset
[-\frac{1}{2},-\frac{1}{16}]\cup[\frac{1}{16},\frac{1}{2}]$ and
$\text{supp }\hat{\psi}_2\subset [-1,1]$ such that
\begin{equation}\label{e:Discrt_Shrlt_Cond_Cone_1}
\sum_{j\geq 0}\abs{\hat{\psi}_1(2^{-2j}\omega)}^2 =1, \;\;\;
\text{for }\abs{\omega}\geq \frac{1}{8}
\end{equation}
and
\begin{equation}\label{e:Discrt_Shrlt_Cond_Cone_2}
\abs{\hat{\psi}_2(\omega-1)}^2+\abs{\hat{\psi}_2(\omega)}^2+\abs{\hat{\psi}_2(\omega+1)}^2=1,
\;\;\; \text{for } \abs{\omega}\leq 1.
\end{equation}
It follows from (\ref{e:Discrt_Shrlt_Cond_Cone_2}) that, for $j\geq
0$,
\begin{equation}\label{e:Discrt_Shrlt_Cond_Cone_3}
\sum_{\ell=-2^j}^{2^j} \abs{\hat{\psi}_2(2^j\omega-\ell)}^2=1,
\;\;\; \text{for }\abs{\omega}\leq 1.
\end{equation}
For a scale index $j\geq 0$ the anisotropic dilation matrices are defined
as
$$A^j_{(1)}= \left(
\begin{array}{cccc}
  4^j    & 0      & \ldots & 0\\
  0      & 2^j    & \ldots & 0\\
  \vdots & \vdots & \vdots & \vdots\\
  0      & 0      & \ldots & 2^j\\
\end{array}
\right), \;\;\; \ldots \;\;\; , A^j_{(d)}=\left(
\begin{array}{cccc}
  2^j    & 0      & \ldots & 0\\
  0      & 2^j      & \ldots & 0\\
  \vdots & \vdots & \vdots & \vdots\\
  0      & 0      & \ldots & 4^j\\
\end{array}
\right),
$$
and for $\ell=(\ell_1,\ldots,\ell_{d-1})$ with
$-2^j\leq\ell_i\leq2^j$, $i=1,\ldots,d-1$, the $d\times d$ shear
matrices are defined as
$$B^{[\ell]}_{(1)}= \left(
\begin{array}{cccc}
  1      & \ell_1 & \ldots & \ell_{d-1}\\
  0      & 1      & \ldots & 0\\
  \vdots & \vdots & \vdots & \vdots\\
  0      & 0      & \ldots & 1\\
\end{array}
\right), \;\;\; \ldots, B^{[\ell]}_{(d)}=\left(
\begin{array}{cccc}
  1      & 0      & \ldots & 0\\
  0      & 1      & \ldots & 0\\
  \vdots & \vdots & \vdots & \vdots\\
\ell_1   & \ell_2 & \ldots & 1\\
\end{array}
\right).
$$
To shorten notation we will write $\abs{[\ell]}\preceq 2^j$ instead
of $\abs{\ell_i}\leq 2^j$, $i=1,\ldots,d-1$, and $\abs{[\ell]}=2^j$
when $\abs{\ell_i}=2^j$ for at least one $i=1,2,\ldots,d-1$. Define
$\hat{\psi}^{(1)}(\xi):=\hat{\psi}_1(\xi_1)\prod_{\frak{d}=2}^{d}\hat{\psi}_2(\frac{\xi_\frak{d}}{\xi_1})$.
Since $\xi A_{(1)}^{-j} B_{(1)}^{[-\ell]}=(4^{-j}\xi_1,
-4^{-j}\xi_1\ell_1+2^{-j}\xi_2, \ldots,
-4^{-j}\xi_1\ell_{d-1}+2^{-j}\xi_d)$, from
(\ref{e:Discrt_Shrlt_Cond_Cone_1}) and
(\ref{e:Discrt_Shrlt_Cond_Cone_3}) it follows that
\begin{eqnarray}\label{e:PrsvlFrm_Prop_Shrlts}
\nonumber
    & &\!\!\!\!\!\!\!\!\!\!\!\!\!\!\!\!\!\!\!\!\!\!\!\!\!\!\!\!\!
        \sum_{j\geq 0}\sum_{\abs{[\ell]}\preceq2^j}
        \abs{\hat{\psi}^{(1)}(\xi A^{-j}_{(1)}B^{[-\ell]}_{(1)})}^2 \\
\nonumber
    &=& \sum_{j\geq0}\sum_{\abs{\ell_1},\ldots,\abs{\ell_{d-1}}\leq2^j}
        \abs{\hat{\psi}_1(2^{-2j}\xi_1)}^2\prod_{\frak{d}=2}^d\abs{\hat{\psi}_2(2^j\frac{\xi_\frak{d}}{\xi_1}-\ell_{\frak{d}-1})}^2 \\
\nonumber
    &=& \sum_{j\geq0} \abs{\hat{\psi}_1(2^{-2j}\xi_1)}^2\sum_{\abs{\ell_1},\ldots,\abs{\ell_{d-2}}\leq2^j}
        \prod_{\frak{d}=2}^{d-1}\abs{\hat{\psi}_2(2^j\frac{\xi_\frak{d}}{\xi_1}-\ell_{\frak{d}-1})}^2\\
\nonumber
    &\vdots&  \\
    &=& 1,
\end{eqnarray}
for $\xi=(\xi_1,\ldots,\xi_d)\in\mathcal{C}^{(1)}$ and which we will
call the \textbf{Parseval frame condition} (for the cone
$\mathcal{C}^{(1)}$). Since $\text{supp }\hat{\psi}^{(1)}\subset
[-\frac{1}{2},\frac{1}{2}]^d$, (\ref{e:PrsvlFrm_Prop_Shrlts})
implies that the shearlet system
\begin{equation}\label{e:Shrlt_Sys_Cone}
\{\psi_{j,\ell,k}^{(1)}(x)= \abs{\text{det
}A_{(1)}}^{j/2}\psi^{(1)}(B^{[\ell]}_{(1)} A^j_{(1)} x-k): j\geq 0,
\abs{[\ell]}\preceq2^j, k\in\mathbb{Z}^d\},
\end{equation}
is a Parseval frame for $L^2((\mathcal{C}^{(1)})^\vee):=\{f\in
L^2(\mathbb{R}^d): \text{supp }\hat{f}\subset\mathcal{C}^{(1)}\}$
(see \cite{GLLWW}, Section 5.2.1). This means that
$$\sum_{j\geq 0}\sum_{\abs{[\ell]}\preceq 2^j}\sum_{k\in\mathbb{Z}^d}
    \abs{\ip{f}{\psi_{j,\ell,k}^{(1)}}}^2 = \norm{f}^2_{L^2(\mathbb{R}^d)},$$
for all $f\in L^2((\mathcal{C}^{(1)})^\vee)$ such that $\text{supp
}\hat{f}\subset \mathcal{C}^{(1)}$. 
One can also construct a shearlet system for any cone
\begin{equation*}
\mathcal{C}^{(i)}=\{\xi=(\xi_1,\ldots,\xi_i,\ldots,\xi_d)\in\hat{\mathbb{R}}^d:
    \abs{\xi_i}\geq\frac{1}{8}, \abs{\frac{\xi_\frak{d}}{\xi_i}}\leq 1, \frak{d}\not=i\},
\end{equation*}
by defining
$\hat{\psi}^{(i)}(\xi)=\hat{\psi}_1(\xi_i)\prod_{\frak{d}\not=i}\hat{\psi}_2(\frac{\xi_\frak{d}}{\xi_i})$
and choosing correspondingly the anisotropic dilation and shear matrices
$A_{(i)}^j$ and $B_{(i)}^{[\ell]}$.

Let $\hat{\Psi}\in C^\infty_c(\mathbb{R}^d)$, with $\text{supp
}\hat{\Psi}\subset [-\frac{1}{4},\frac{1}{4}]^d$ and
$\abs{\hat{\Psi}}=1$ for $\xi\in
[-\frac{1}{8},\frac{1}{8}]^d=\mathcal{R}$, be such that
\begin{equation}\label{e:ShrltPrsvlFrm}
   \abs{\hat{\Psi}(\xi)}^2\chi_\mathcal{R}(\xi)
        + \sum_{\frak{d}=1}^d\sum_{j\geq 0}\sum_{\abs{[\ell]}\preceq2^j}\abs{\hat{\psi}^{(\frak{d})}(\xi A^{-j}_{(\frak{d})}B^{[-\ell]}_{(\frak{d})})}^2
        \chi_{\mathcal{C}^{(\mathfrak{d})}}(\xi)=1,
\end{equation}
for all $\xi\in\hat{\mathbb{R}}^d$. This implies that one can
construct a Parseval frame for $L^2(\mathbb{R}^d)$ (see Theorem 9 in \cite{LaWe09}).

Since $\mathcal{C}^{(i)}$ are orthogonal rotations of
$\mathcal{C}^{(1)}$ we will often drop the sub- or super- index and
develop our results only for one direction. We will incorporate the
directions only in the definitions of the spaces. 

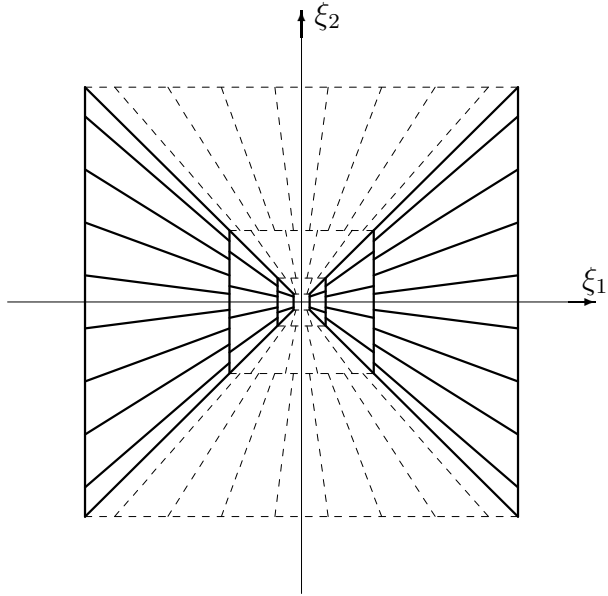
\begin{figure}
\vskip0.8cm
\begin{center}
\unitlength=1pt
\begin{picture}(220,220)

\drawline(110,0)(110,220)\put(110,210){\vector(0,1){10}}

\put(115,215){$\xi_2$}

\drawline(0,110)(220,110)\put(210,110){\vector(1,0){10}}

\put(215,115){$\xi_1$}

\thicklines



\drawline(113,107)(113,113) \drawline(119,101)(119,119)
\drawline(137, 83)(137,137) \drawline(191, 29)(191, 191)

\drawline(113,113)(191,191) \drawline(113,107)(191,29)

\drawline(113,112)(119,114) \drawline(113,108)(119,106)

\drawline(119,116)(137,129) \drawline(119,112)(137,116)
\drawline(119,108)(137,104) \drawline(119,104)(137,91)

\drawline(137,133)(191,180) \drawline(137,126)(191,160)
\drawline(137,120)(191,140) \drawline(137,113)(191,120)
\drawline(137,107)(191,100) \drawline(137,100)(191,80)
\drawline(137,94)(191,60) \drawline(137,87)(191,40)


\drawline(107,107)(107,113) \drawline(101,101)(101,119)
\drawline(83, 83)(83,137) \drawline(29, 29)(29, 191)

\drawline(29,191)(107,113) \drawline(29,29)(107,107)

\drawline(107,112)(101,114) \drawline(107,108)(101,106)

\drawline(101,116)(83,129) \drawline(101,112)(83,116)
\drawline(101,108)(83,104) \drawline(101,104)(83,91)

\drawline(83,133)(29,180) \drawline(83,126)(29,160)
\drawline(83,120)(29,140) \drawline(83,113)(29,120)
\drawline(83,107)(29,100) \drawline(83,100)(29,80)
\drawline(83,94)(29,60) \drawline(83,87)(29,40)



\thinlines

\dashline[+30]{3}(107,113)(113,113)
\dashline[+30]{3}(101,119)(119,119)
\dashline[+30]{3}(83,137)(137,137)
\dashline[+30]{3}(29,191)(191,191)


\dashline[+30]{3}(112,113)(114,119)
\dashline[+30]{3}(108,113)(106,119)

\dashline[+30]{3}(116,119)(129,137)
\dashline[+30]{3}(112,119)(116,137)
\dashline[+30]{3}(108,119)(104,137)
\dashline[+30]{3}(104,119)(91,137)

\dashline[+30]{3}(133,137)(180,191)
\dashline[+30]{3}(126,137)(160,191)
\dashline[+30]{3}(120,137)(140,191)
\dashline[+30]{3}(113,137)(120,191)
\dashline[+30]{3}(107,137)(100,191)
\dashline[+30]{3}(100,137)(80,191) \dashline[+30]{3}(94,137)(60,191)
\dashline[+30]{3}(87,137)(40,191)

%
\dashline[+30]{3}(107,107)(113,107)
\dashline[+30]{3}(101,101)(119,101) \dashline[+30]{3}(83,83)(137,83)
\dashline[+30]{3}(29,29)(191,29)


\dashline[+30]{3}(112,107)(114,101)
\dashline[+30]{3}(108,107)(106,101)

\dashline[+30]{3}(116,101)(129,83)
\dashline[+30]{3}(112,101)(116,83)
\dashline[+30]{3}(108,101)(104,83) \dashline[+30]{3}(104,101)(91,83)

\dashline[+30]{3}(133,83)(180,29) \dashline[+30]{3}(126,83)(160,29)
\dashline[+30]{3}(120,83)(140,29) \dashline[+30]{3}(113,83)(120,29)
\dashline[+30]{3}(107,83)(100,29) \dashline[+30]{3}(100,83)(80,29)
\dashline[+30]{3}(94,83)(60,29) \dashline[+30]{3}(87,83)(40,29)

\end{picture}
\end{center}
\caption{Sketch of the partition of the frequency plane
$\hat{\mathbb{R}}^2$ induced by the
shearlets.}\label{f:tiling_freq_shear} \vskip0.8cm
\end{figure}

Observe that the characteristic functions $\chi_R$ and
$\chi_{\mathcal{D}^{(\mathfrak{d})}}$ in (\ref{e:ShrltPrsvlFrm})
destroy the localization in the space domain. The characteristic functions $\chi_R$
and $\chi_{\mathcal{D}^{(\mathfrak{d})}}$ can be removed from (\ref{e:ShrltPrsvlFrm})
in whose case the condition of tight frame with bounds equal 1 (Parseval frame) will be lost,
however the property of being a frame will remain. With a slight
variation on the above, Guo and Labate constructed in \cite{GL11}
smooth Parseval frames of shearlets on the cone (whose definition is
not important in the development of this paper), which means that a)
one can ignore the characteristic functions $\chi_R,
\chi_{\mathcal{D}^{(\mathfrak{d})}}$ in (\ref{e:ShrltPrsvlFrm}) and
b) one can prove that the composition of the analysis (see
(\ref{e:AnlsOpShrlt})) and synthesis operators is the identity on
$\mathcal{S}'$ (see \cite{Ver12}). The use of the frame of shearlets
(without the characteristic functions $\chi_R$ and
$\chi_{\mathcal{D}^{(\mathfrak{d})}}$)
or the smooth Parseval frame of shearlets is
transparent to the results in this paper.

Since $\mathcal{D}^{(\mathfrak{d})}$ are orthogonal rotations of
$\mathcal{D}^{(1)}$ we will often drop the sub- or super- index and
develop our results only for one direction. 

\begin{figure}
\vskip0.8cm
\begin{center}
\unitlength=0.7pt
\begin{picture}(360,220)

\drawline(0,20)(360,20)\put(350,20){\vector(1,0){10}}

\put(350,10){$x_1$}

\drawline(20,0)(20,220)\put(20,210){\vector(0,1){10}}

\put(25,210){$x_2$}

\dashline[+30]{3}(0,180)(360,180) \dashline[+30]{3}(0,100)(360,100)

\dashline[+30]{3}(0,200)(20,220)\dashline[+30]{3}(0,160)(60,220)
\dashline[+30]{3}(0,120)(100,220)\dashline[+30]{3}(0,80)(140,220)
\dashline[+30]{3}(0,40)(180,220) \dashline[+30]{3}(0,0)(220,220)
\dashline[+30]{3}(40,0)(260,220) \dashline[+30]{3}(80,0)(300,220)
\dashline[+30]{3}(120,0)(340,220)
\dashline[+30]{3}(160,0)(360,200)
\dashline[+30]{3}(200,0)(360,160)\dashline[+30]{3}(240,0)(360,120)
\dashline[+30]{3}(280,0)(360,80)\dashline[+30]{3}(320,0)(360,40)

\thicklines

\drawline(20,20)(180,180) \drawline(20,20)(180,20)
\drawline(180,20)(340,180)\drawline(180,180)(340,180)

\drawline(60,20)(220,180) \drawline(100,20)(260,180)
\drawline(140,20)(300,180)

\drawline(100,100)(260,100)

\put(18,180){\line(1,0){4}}\put(12,185){$1$}
\put(180,18){\line(0,1){4}}\put(182,10){$1$}
\put(340,18){\line(0,1){4}}\put(342,10){$2$}

\end{picture}
\end{center}
\caption{Sketch of the covering of the plane $\mathbb{R}^2$ with
parallelograms in $\mathcal{Q}^{1,2}$. Those parallelograms with
solid lines cover the parallelogram
$Q_{j,\ell,k}=Q_{0,1,0}$.}\label{f:Covering_Real_Plane} 
\end{figure}
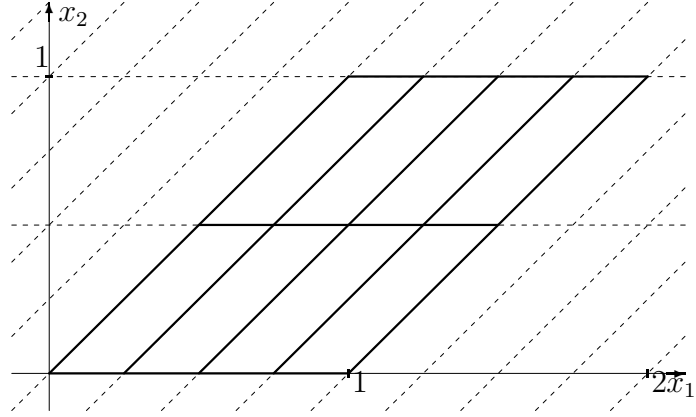

For $Q_0=[0,1)^d$, write
\begin{equation}\label{e:def_Qjlk}
Q^{(\mathfrak{d})}_{j,\ell,k}=A^{-j}_{(\mathfrak{d})}B_{(\mathfrak{d})}^{[-\ell]}(Q_0+k),
\end{equation}
with $\mathfrak{d}=1,\ldots,d$, $j\geq 0$, $\abs{[\ell]}\preceq 2^j$
and $k\in\mathbb{Z}^d$. Therefore, $\int
\chi_{Q_{j,\ell,k}^{(\mathfrak{d})}}=\abs{Q^{(\mathfrak{d})}_{j,\ell,k}}=2^{-(d+1)j}=\abs{\text{det
} A_{(\mathfrak{d})}}^{-j}$. Let
$\mathcal{Q}_{AB}:=\{Q^{(\mathfrak{d})}_{j,\ell,k}:
\mathfrak{d}=1,\ldots,d, j\geq 0,\abs{[\ell]}\preceq 2^j,
k\in\mathbb{Z}^d\}$ and $\mathcal{Q}^{j,\ell}_{(\mathfrak{d})}:=
\{Q^{(\mathfrak{d})}_{j,\ell,k}: k\in\mathbb{Z}^d\}$. Thus, for
fixed $\mathfrak{d}, j, \ell$,
$\mathcal{Q}^{j,\ell}_{(\mathfrak{d})}$ is a partition of
$\mathbb{R}^d$ as can be seen in Figure \ref{f:Covering_Real_Plane}.
Hence, for every $j\geq0$ there exist $2^{(j+1)(d-1)}+1$ partitions
since $\abs{[\ell]}\preceq2^j$. To shorten notation and clear
exposition, we will identify the multi indices $(j,\ell,k)$ and
$(i,m,n)$ with $P$ and $Q$, respectively. This way we write
$\psi_P=\psi_{j,\ell,k}$ or $\psi_Q=\psi_{i,m,n}$, regardless the
direction in question. We also write
$\tilde{\chi}_Q(x)=\abs{Q}^{-1/2}\chi_Q(x)$. 
We formally define the shearlet \textbf{analysis operator} as
\begin{equation}\label{e:AnlsOpShrlt}
  S_{\Psi,\psi} f = \{ \{\ip{f}{\Psi(\cdot - k)}\}_{k\in\mathbb{Z}^d} ,
  \{\ip{f}{\psi_Q}\}_{Q\in\mathcal{Q}_{AB}}\}.
\end{equation}
So the shearlet coefficients are $\mathbf{s}=
\{s_Q\}_{Q\in\mathcal{Q}_{AB}}= S_{\Psi,\psi} f$.

The next two definitions are the sequence spaces associated to the
shear anisotropic inhomogeneous Besov and Triebel-Lizorkin
distribution spaces as defined in \cite{Ver13} and \cite{Ver12},
respectively.

For $s\in \mathbb{R}$, $0<p,q\leq\infty$, the \textbf{shear
anisotropic inhomogeneous Besov sequence space}
$\mathbf{b}^{s,q}_p(AB)$ is defined as the collection of all
complex-valued sequences $\mathbf{s}=\{s_Q\}_{Q\in\mathcal{Q}_{AB}}$
such that
\begin{eqnarray}\label{e:def_BsvSeq_Shrlts}
\nonumber
 \norm{\mathbf{s}}_{\mathbf{b}^{s,q}_p(AB)}
    &:=& \left(\sum_{k\in\mathbb{Z}^d} \abs{s_k}^p\right)^{1/p}\\
    &+& \left(\sum_{\mathfrak{d}=1}^d\sum_{j\geq 0}\sum_{\abs{[\ell]}\preceq
        2^j}\left(\sum_{Q\in\mathcal{Q}^{j,\ell}_{(\mathfrak{d})}}
        [\abs{Q}^{-s+\frac{1}{p}-\frac{1}{2}}\abs{s_Q}]^p\right)^{q/p}\right)^{1/q}<\infty.
\end{eqnarray}

For $s\in\mathbb{R}$, $0<p<\infty$ and $0<q\leq \infty$.  The
\textbf{shear anisotropic inhomogeneous Triebel-Lizorkin sequence
space} $\mathbf{f}^{s,q}_p(AB)$ is defined as the collection of all
complex-valued sequences $\mathbf{s}=\{s_Q\}_{Q\in\mathcal{Q}_{AB}}$
such that
\begin{equation}\label{e:def_T-L_seq_spcs}
  \norm{\mathbf{s}}_{\mathbf{f}^{s,q}_p(AB)}
    = \left(\sum_{k\in\mathbb{Z}^d} \abs{s_k}^p\right)^{1/p} +
    \norm{\left(\sum_{Q\in\mathcal{Q}_{AB}}(\abs{Q}^{-s}\abs{s_Q}\tilde{\chi}_Q)^q\right)^{1/q}}_{L^p}<\infty.
\end{equation}

It is not hard to prove that, for $0<p<\infty$,
$\mathbf{f}^{s,p}_p(AB)=\mathbf{b}^{s,p}_p(AB)$, either by straight
calculations from the definitions above or by the embedding Theorem
4.1 iii) in \cite{Ver13} with $p=q$.

Next result is a particular case which identifies the shear
anisotropic Triebel-Lizorkin or Besov sequence spaces just defined
with discrete Lorentz sequence spaces indexed by $\mathcal{Q}_{AB}$
for certain parameters.

\begin{lem}\label{l:BsvSeqSpcs_r_LrntzSpcs}
For $\beta, s\in\mathbb{R}$ and $0<\tau<\infty$, let
$\gamma=s+\frac{1-\beta}{\tau}$ and
$\mathbf{u}=\{\abs{Q}^{-s-\frac{1}{2}}\}_{Q\in\mathcal{Q}_{AB}}$.
Then,
$$\ell^{\tau,\tau}(\mathbf{u},\nu_\beta, \mathcal{Q}_{AB})=\mathbf{b}^{\gamma,\tau}_\tau(AB)=\mathbf{f}^{\gamma,\tau}_\tau(AB),$$
with equal quasi-norms.
\end{lem}
\textbf{Proof}. This is a straight consequence of the respective
definitions and the conditions on the parameters. So,
\begin{eqnarray*}
  \norm{\mathbf{s}}_{\ell^{\tau,\tau}(\mathbf{u},\nu_\beta,\mathcal{Q}_{AB})}^\tau
    &=& \norm{\{\abs{s_Q} \abs{Q}^{-s-\frac{1}{2}}\}_{Q\in\mathcal{Q}_{AB}}}_{\ell^{\tau,\tau}(\nu_\beta,\mathcal{Q}_{AB})}^\tau \\
    &=& \sum_{Q\in\mathcal{Q}_{AB}} (\abs{s_Q}\abs{Q}^{-s-\frac{1}{2}})^\tau
        \abs{Q}^\beta
        = \sum_{Q\in\mathcal{Q}_{AB}}
        (\abs{s_Q}\abs{Q}^{-\gamma+\frac{1}{\tau}-\frac{1}{2}})^\tau \\
    &=& \norm{\mathbf{s}}_{\mathbf{b}^{\gamma,\tau}_\tau(AB)}^\tau.
\end{eqnarray*}
\hfill $\blacksquare$ \vskip .5cm   

\vskip 1cm
\section{Democracy of $\mathbf{b}^{s,p}_p(AB)$ and $\mathbf{f}^{s,q}_p(AB)$.}\label{S:Upr-Lwr_Tmlkv-prop_T-L_Shrlts}
Our aim now is to prove that the spaces $\mathbf{b}^{s,p}_p(AB)$ and
$\mathbf{f}^{s,q}_p(AB)$ verify points J1) and B1) of Theorem
\ref{t:J-B_ineq_equiv_for_seq}, with $\mathbf{u}$ related to a
second space. For $\mathbf{f}^{s,q}_p(AB)$ we will not have the same
upper and lower democracy. Hence, only embeddings can be proved.
However, for $\mathbf{f}^{s,p}_p(AB)=\mathbf{b}^{s,p}_p(AB)$ the
upper and lower democracy are the same and this fact allows us to
prove full characterizations.


To prove democracy for $\mathbf{f}^{s,q}_p(AB)$ we need a previous
result.
\begin{lem}\label{l:S_gamma_bnd_T-L_Shrlts}
Let $S^\gamma_\Gamma(x):=\sum_{P\in\Gamma}\abs{P}^\gamma\chi_P(x)$,
$P^x$ be the largest ``cube" that contains $x$ for some scale and
$P_x$ be a smallest ``cube" that contains $x$ for some scale.
\begin{itemize}
  \item[a)] If $\gamma>\frac{d-1}{d+1}$ and there exists $P^x$ then,
        $$\abs{P^x}^\gamma\chi_{P^x}(x)\leq S^\gamma_\Gamma(x)\leq
        C_\gamma\abs{P^x}^{\gamma-\frac{d-1}{d+1}}\chi_{P^x(x)}.$$
  \item[b)] If $\gamma<\frac{d-1}{d+1}$ and there exists $P_x$ then,
        $$\abs{P_x}^\gamma\chi_{P_x}(x)\leq S^\gamma_\Gamma(x)\leq
        C_\gamma\abs{P_x}^{\gamma-\frac{d-1}{d+1}}\chi_{P_x(x)}.$$
\end{itemize}
\end{lem}
\begin{rem}\label{r:l:S_gamma_bnd_T-L_Shrlts}
Observe that we refer to \textbf{the} largest ``cube" $P^x$ and to
\textbf{a} smallest ``cube" $P_x$. For any scale $j\geq0$ there are
$2^{(j+1)(d-1)}+1$ different coverings of the plane $\mathbb{R}^d$ since
$\abs{[\ell]}\preceq 2^j$. Therefore, a ``child" cube that contain
$x$ will not be unique.
\end{rem}
\textbf{Proof of Lemma \ref{l:S_gamma_bnd_T-L_Shrlts}}. We start by
proving a). It is clear that $\abs{P^x}^\gamma\chi_{P^x}(x)\leq
S^\gamma_\Gamma(x)$, since the sum in $S^\gamma_\Gamma(x)$ contains
$\abs{P^x}\chi_{P^x}(x)$, at least. For the right-hand side of the
inequality we enlarge the sum defining $S^\gamma_\Gamma(x)$ to
include all $P$'s in the same and finer scales that contain $x$.
Defining $\Gamma^{j,\ell}:=\Gamma\cap\mathcal{Q}^{j,\ell}$ and since
$\mathcal{Q}^{j,\ell}$ is a partition of $\mathbb{R}^2$, we obtain
\begin{eqnarray*}
  S^\gamma_\Gamma(x)
    &=& \sum_{j\geq J}\sum_{\abs{[\ell]}\preceq2^j}\sum_{P\in\Gamma^{j,\ell}}\abs{P}^\gamma\chi_P(x)
        \leq \sum_{j\geq J}\sum_{\abs{[\ell]}\preceq2^j}\sum_{P\in\mathcal{Q}^{j,\ell}}\abs{P}^\gamma\chi_P(x) \\
    &=& \sum_{j\geq J} \abs{P_j}^\gamma (2^{(j+1)(d-1)}+1)\chi_{P^x}(x)
        \leq C_d \sum_{j\geq J} 2^{-j((d+1)\gamma-(d-1))}\chi_{P^x}(x) \\
    &=& C_{d,\gamma}\abs{P^x}^{\gamma-\frac{(d-1)}{(d+1)}}\chi_{P^x}(x),
\end{eqnarray*}
since $\gamma>\frac{d-1}{d+1}$ and $\abs{P^x}=2^{-J(d+1)}$.

The left-hand side of b) is clear since $\Gamma\ni P_x$, at least.
For the right-hand side of b) we enlarge the sum
$S^\gamma_\Gamma(x)$ to include all ``cubes" in the same and coarser
scales. With the same definition of $\Gamma^{j,\ell}$ and since
$\mathcal{Q}^{j,\ell}$ is a partition of $\mathbb{R}^d$, for fixed
$j,\ell$, we obtain
\begin{eqnarray*}
  S^\gamma_\Gamma(x)
    &=& \sum_{j=0}^J\sum_{\abs{[\ell]}\preceq 2^j}\sum_{P\in\Gamma^{j,\ell}} \abs{P}^\gamma\chi_P(x)
        \leq \sum_{j=0}^J\sum_{\abs{[\ell]}\preceq 2^j}\sum_{P\in\mathcal{Q}^{j,\ell}} \abs{P}^\gamma\chi_P(x) \\
    &=& \sum_{j=0}^J \abs{P_j}^\gamma (2^{(j+1)(d-1)}+1)\chi_{P_x}(x)
        \leq C_d \sum_{j=0}^J 2^{-j[(d+1)\gamma-(d-1)]} \chi_{P_x}(x) \\
    &\leq& C_d 2^{-J[(d+1)\gamma-(d-1)]} \sum_{j=0}^\infty 2^{+j[(d+1)\gamma-(d-1)]}\chi_{P_x}(x)
        = C_{d,\gamma} \abs{P_x}^{\gamma-\frac{d-1}{d+1}}\chi_{P_x}(x),
\end{eqnarray*}
because $\gamma<\frac{d-1}{d+1}$.

\hfill $\blacksquare$ \vskip .5cm   

The first main result is:
\begin{thm}\label{t:Tmlkv_T-L_Shrlts}
Let $s_1,s_2\in \mathbb{R}$, $0<p_1,p_2<\infty$,
$0<q_1,q_2\leq\infty$ with $p_1\neq q_1$ and
$\mathbf{u}=\{\norm{\mathbf{e}_P}_{\mathbf{f}_{p_2}^{s_2,q_2}(AB)}\}_{P\in\mathcal{Q}_{AB}}=\{\abs{P}^{-s_2+\frac{1}{p_2}-\frac{1}{2}}\}_{P\in\mathcal{Q}_{AB}}$.
If $\alpha=p_1(s_2-\frac{1}{p_2}-s_1+\frac{1}{p_1})$, there exist
$C,C'>0$ depending only on $s_1,s_2,p_2$ and $q_1$ such that
\begin{equation}\label{e:Tmlkv_T-L_Shrlts}
C(\nu_{\alpha+\frac{d-1}{d+1}}(\Gamma))^{1/p_1} \leq
    \norm{\sum_{P\in\Gamma} \frac{\mathbf{e}_P}{\norm{\mathbf{e}_P}_{\mathbf{f}_{p_2}^{s_2,q_2}(AB)}}}_{\mathbf{f}_{p_1}^{s_1,q_1}(AB)}
    \leq
    C'(\nu_{\alpha-\frac{p_1(d-1)}{q_1(d+1)}}(\Gamma))^{1/p_1}
\end{equation}
for all $\Gamma\subset\mathcal{Q}_{AB}$ such that
$\nu_{\alpha-\frac{p_1(d-1)}{q_1(d+1)}}(\Gamma)<\infty$. Conversely,
if (\ref{e:Tmlkv_T-L_Shrlts}) holds then,
$$p_1(s_2-\frac{1}{p_2}-s_1+\frac{1}{p_1})-\frac{p_1(d-1)}{q_1(d+1)}
    \leq \alpha\leq
    p_1(s_2-\frac{1}{p_2}-s_1+\frac{1}{p_1}) +
    \frac{d-1}{d+1}.$$
\end{thm}

\textbf{Proof}. Write $\mathbf{f}_1:=\mathbf{f}_{p_1}^{s_1,q_1}(AB)$
and $\mathbf{f}_2:=\mathbf{f}_{p_2}^{s_2,q_2}(AB)$ to simplify the
notation in this proof.  By definition we have that
\begin{eqnarray*}
  \norm{\sum_{P\in\Gamma}\frac{\mathbf{e}_P}{\norm{\mathbf{e}_P}_{\mathbf{f}_2}}}_{\mathbf{f}_1}
    &=& \left(\int_{\mathbb{R}^2} \left[
        \sum_{j\geq 0}\sum_{\abs{[\ell]}\preceq 2^j}\sum_{P\in \mathcal{Q}^{j,\ell}}
            (\abs{P}^{s_2-s_1-\frac{1}{p_2}}\chi_{P}(x) )^{q_1}\right]^{p_1/q_1}
            dx\right)^{1/p_1}.
\end{eqnarray*}

We consider three cases:

1) We start by assuming that
$\alpha=p_1(s_2-\frac{1}{p_2}-s_1+\frac{1}{p_1})>1+\frac{p_1(d-1)}{q_1(d+1)}
\Leftrightarrow \gamma=q_1(s_2-s_1-\frac{1}{p_2})>\frac{d-1}{d+1}$.
In this case, since
$\nu_{\alpha-\frac{p_1(d-1)}{q_1(d+1)}}(\Gamma)<\infty$, the largest
$P^x\in\Gamma$ exists for all $x\in\cup_{P\in\Gamma}P$. Applying
part a) of Lemma \ref{l:S_gamma_bnd_T-L_Shrlts}, we get
\begin{eqnarray*}
    && \!\!\!\!\!\!\!\!\!\!\!\!\!\!\!\!\!\!\!\!\!\!\!\!
            \left[\sum_{j\geq 0}\sum_{\abs{[\ell]}\preceq 2^j}\sum_{P\in\Gamma^{j,\ell}}
            \abs{P}^{q_1(s_2-s_1-\frac{1}{p_2})} \chi_P(x)\right]^{p_1/q_1} \\
    &\leq& \left[C_\gamma \abs{P^x}^{\gamma-\frac{d-1}{d+1}}\chi_{P^x}(x)\right]^{p_1/q_1} \\
    &=& C_{\gamma,p_1} \abs{P^x}^{\frac{p_1}{q_1}(\gamma-\frac{d-1}{d+1})}\chi_{P^x}(x) \\
    &\leq& C_{\gamma,p_1} \sum_{j\geq 0}\sum_{\abs{[\ell]}\preceq
        2^j}\sum_{P\in\Gamma^{j,\ell}} \abs{P}^{\frac{p_1}{q_1}(\gamma-\frac{d-1}{d+1})}\chi_{P}(x)\\
    &=& C_{\gamma,p_1} \sum_{P\in\Gamma}
    \abs{P}^{\frac{p_1}{q_1}(\gamma-\frac{d-1}{d+1})}\chi_P(x),
\end{eqnarray*}
for all $x\in\cup_{P\in\Gamma}P$. From this we deduce the
upper-Temlyakov property as
\begin{eqnarray*}
  \norm{\sum_{P\in\Gamma}\frac{\mathbf{e}_P}{\norm{\mathbf{e}_P}_{\mathbf{f}_2}}}_{\mathbf{f}_1}
    &\leq& C_\gamma \left(\int_{\mathbb{R}^2} \sum_{P\in\Gamma} \abs{P}^{p_1(s_2-s_1-\frac{1}{p_2})-\frac{p_1(d-1)}{q_1(d+1)}} \chi_P(x) dx\right)^{1/p_1} \\
    &=& C_\gamma \left(\sum_{P\in\Gamma} \abs{P}^{p_1(s_2-s_1-\frac{1}{p_2})+1-\frac{p_1(d-1)}{q_1(d+1)}} \right)^{1/p_1} \\
    &=& C_\gamma(\nu_{\alpha-\frac{p_1(d-1)}{q_1(d+1)}}(\Gamma))^{1/p_1}.
\end{eqnarray*}
We also have from part a) of Lemma \ref{l:S_gamma_bnd_T-L_Shrlts}
that
\begin{eqnarray*}
    & & \!\!\!\!\!\!\!\!\!\!\!\!\!\!\!\!\!\!\!\!\!\!\!\!
            \left[\sum_{j\geq 0}\sum_{\abs{\ell}\leq 2^j}\sum_{P\in\Gamma^{j,\ell}}
            \abs{P}^{q_1(s_2-s_1-\frac{1}{p_2})}
            \chi_P(x)\right]^{p_1/q_1} \\
    &\geq&  \left[\abs{P^x}^{q_1(s_2-s_1-\frac{1}{p_2})} \chi_{P^x}(x)\right]^{p_1/q_1} \\
    &=& \abs{P^x}^{p_1(s_2-s_1-\frac{1}{p_2})} \chi_{P^x}(x) \\
    &\geq& \frac{1}{C_\gamma} \sum_{P\in\Gamma}
    \abs{P}^{p_1(s_2-s_1-\frac{1}{p_2})+\frac{d-1}{d+1}}\chi_P(x).
\end{eqnarray*}
From this we obtain the lower-Temlyakov property as
\begin{eqnarray*}
  \norm{\sum_{P\in\Gamma}\frac{\mathbf{e}_P}{\norm{\mathbf{e}_P}_{\mathbf{f}_2}}}_{\mathbf{f}_1}
    &\geq& \left(\frac{1}{C_\gamma}\int_{\mathbb{R}^2} \sum_{P\in\Gamma} \abs{P}^{p_1(s_2-s_1-\frac{1}{p_2})+\frac{d-1}{d+1}}\chi_P(x) dx \right)^{1/p_1} \\
    &=& \frac{1}{C_{\gamma,p_1}}\left(\sum_{P\in\Gamma} \abs{P}^{p_1(s_2-s_1-\frac{1}{p_2})+1+\frac{d-1}{d+1}} \right)^{1/p_1} \\
    &=& \frac{1}{C_{\gamma,p_1}}\left(\sum_{P\in\Gamma} \abs{P}^{\alpha+\frac{d-1}{d+1}}
    \right)^{1/p_1}.
\end{eqnarray*}

2) Consider now the case
$\alpha=p_1(s_2-s_1-\frac{1}{p_2})+1<1+\frac{p_1(d-1)}{q_1(d+1)}
\Leftrightarrow \gamma=q_1(s_2-s_1-\frac{1}{p_2})<\frac{d-1}{d+1}$.
We can show that the set $E_\alpha$ of all $x\in\cup_{P\in\Gamma}P$
for which $P_x$ does not exist has measure zero. To see this, we
remind that $\mathcal{Q}^{j,\ell}=\{P\in\mathcal{Q}_{AB}:
\abs{P}=\abs{Q_{j,\ell,k}}=\abs{Q_j}=2^{-j(d+1)}, j\geq 0\}$. Then,
for all $i\geq 0$, $E_\alpha\subset \cup_{j\geq
i}\cup_{\abs{[\ell]}\preceq 2^j}\cup_{P\in\Gamma^{j,\ell}}P$.
Therefore,
\begin{eqnarray*}
  \abs{E_\alpha}
    &\leq& \sum_{j\geq i}\sum_{\abs{[\ell]}\preceq 2^j}\sum_{P\in\Gamma^{j,\ell}} \abs{P} \\
    &=& \sum_{j\geq i}\sum_{\abs{[\ell]}\preceq 2^j}\sum_{P\in\Gamma^{j,\ell}} \abs{P}^{\alpha-\frac{p_1(d-1)}{q_1(d+1)}}\abs{P}^{1-(\alpha-\frac{p_1(d-1)}{q_1(d+1)})} \\
    &\leq& \sum_{j\geq i} 2^{-j(d+1)(1+\frac{p_1(d-1)}{q_1(d+1)}-\alpha)} \nu_{\alpha-\frac{p_1(d-1)}{q_1(d+1)}}(\Gamma^j) \\
    &\leq& C_{d,p_1,q_1,\alpha}\nu_{\alpha-\frac{p_1(d-1)}{q_1(d+1)}}(\Gamma)
    2^{-i(d+1)(1+\frac{p_1(d-1)}{q_1(d+1)}-\alpha)},
\end{eqnarray*}
since $\alpha<1+\frac{p_1(d-1)}{q_1(d+1)}$, letting
$i\rightarrow\infty$ we deduce $\abs{E_\alpha}=0$. We now apply
Lemma \ref{l:S_gamma_bnd_T-L_Shrlts}, part b), to obtain
\begin{eqnarray*}
    & & \!\!\!\!\!\!\!\!\!\!\!\!\!\!\!\!\!\!\!\!
        \left[\sum_{j\geq 0}\sum_{\abs{[\ell]}\preceq 2^j}\sum_{P\in\Gamma^{j,\ell}}
        \abs{P}^{q_1(s_2-s_1-\frac{1}{p_2})}\chi_P(x)\right]^{p_1/q_1}\\
    &\leq& C_\gamma\abs{P_x}^{(\gamma-\frac{d-1}{d+1})\frac{p_1}{q_1}}\chi_{P_x}(x) \\
    &\leq& C_\gamma \sum_{j\geq 0}\sum_{\abs{[\ell]}\preceq 2^j}\sum_{P\in\Gamma^{j,\ell}}
        \abs{P}^{p_1(s_2-s_1-\frac{1}{p_2})-\frac{p_1(d-1)}{q_1(d+1)}}\chi_P(x),
\end{eqnarray*}
for all $x\in\cup_{P\in\Gamma}P\setminus E_\alpha$. From this we
deduce the upper-Temlyakov property as
\begin{eqnarray*}
  \norm{\sum_{P\in\Gamma}\frac{\mathbf{e}_P}{\norm{\mathbf{e}_P}_{\mathbf{f}_2}}}_{\mathbf{f}_1}
    &\leq& C_\gamma\left(\sum_{P\in\Gamma}\abs{P}^{p_1(s_2-s_1-\frac{1}{p_2})+1-\frac{p_1(d-1)}{q_1(d+1)}}\right)^{1/p_1} \\
    &=& C_\gamma\left(\sum_{P\in\Gamma}\abs{P}^{\alpha-\frac{p_1(d-1)}{q_1(d+1)}}\right)^{1/p_1}
        = C_\gamma \left(\nu_{\alpha-\frac{p_1(d-1)}{q_1(d+1)}}(\Gamma)\right)^{1/p_1}.
\end{eqnarray*}
For the lower-Temlyakov property we consider two cases:
$\gamma\leq0$ and $0<\gamma<\frac{d-1}{d+1}$. From
$\gamma=q_1(s_2-s_1-\frac{1}{p_2})>0$, it follows that
$s_2-s_1-\frac{1}{p_2}>0\Rightarrow
p_1(s_2-s_1-\frac{1}{p_2})+\frac{d-1}{d+1}>\frac{d-1}{d+1}$. Hence,
part a) of Lemma \ref{l:S_gamma_bnd_T-L_Shrlts} yields
\begin{eqnarray*}
    & & \!\!\!\!\!\!\!\!\!\!\!\!\!\!\!\!\!\!\!\!\!\!\!\!\!\!\!\!\!\!\!\!\!\!\!\!\!
        \left[\sum_{j\geq 0}\sum_{\abs{[\ell]}\preceq 2^j}\sum_{P\in\Gamma^{j,\ell}}
            \abs{P}^{q_1(s_2-s_1-\frac{1}{p_2})}\chi_P(x)\right]^{\frac{p_1}{q_1}} \\
    &\geq& \abs{P^x}^{p_1(s_2-s_1-\frac{1}{p_2})}\chi_{P^x}(x) \\
    &\geq& \frac{1}{C_\gamma}\sum_{P\in\Gamma}\abs{P}^{p_1(s_2-s_1-\frac{1}{p_2})+\frac{d-1}{d+1}}\chi_P(x).
\end{eqnarray*}
Similarly for $\gamma=q_1(s_2-s_1-\frac{1}{p_2})\leq0$,
$s_2-s_1-\frac{1}{p_2}\leq0\Rightarrow
p_1(s_2-s_1-\frac{1}{p_2})+\frac{d-1}{d+1}\leq\frac{d-1}{d+1}$.
Hence, part b) of Lemma \ref{l:S_gamma_bnd_T-L_Shrlts} yields the
same lower bound. In both cases we obtain
\begin{eqnarray*}
  \norm{\sum_{P\in\Gamma}\frac{\mathbf{e}_P}{\norm{\mathbf{e}_P}_{\mathbf{f}_2}}}_{\mathbf{f}_1}
    &\geq& \left(\frac{1}{C_\gamma}\int_{\mathbb{R}^d} \sum_{P\in\Gamma}
        \abs{P}^{p_1(s_2-s_1-\frac{1}{p_2})+\frac{d-1}{d+1}}\chi_P(x)dx\right)^{1/p_1} \\
    &=& \frac{1}{C_{\gamma}}\left(\sum_{P\in\Gamma}\abs{P}^{p_1(s_2-s_1-\frac{1}{p_2})+1+\frac{d-1}{d+1}}\right)^{1/p_1} \\
    &=& \frac{1}{C_{\gamma}}\left(\sum_{P\in\Gamma}\abs{P}^{\alpha+\frac{d-1}{d+1}}\right)^{1/p_1}
        = \frac{1}{C_{\gamma}}\left(\nu_{\alpha+\frac{d-1}{d+1}}(\Gamma)\right)^{1/p_1}.
\end{eqnarray*}

3) For
$\alpha=p_1(s_2-\frac{1}{p_2}-s_1+\frac{1}{p_1})=1+\frac{p_1(d-1)}{q_1(d+1)}
\Leftrightarrow \gamma=q_1(s_2-s_1-\frac{1}{p_2})=\frac{d-1}{d+1}$,
the set $E_{\alpha-\frac{p_1(d-1)}{q_1(d+1)}}$ of all
$x\in\cup_{P\in\Gamma}Q$ for which $P_x$ does not exist has also
measure zero. Indeed, since $\alpha-\frac{p_1(d-1)}{q_1(d+1)}=1$,
$$\abs{E_{\alpha-\frac{p_1(d-1)}{q_1(d+1)}}}\leq
    \sum_{j\geq i}\sum_{\abs{[\ell]}\preceq2^j}\sum_{P\in\Gamma^{j,\ell}}\abs{P}
    = \sum_{j\geq i}\nu_1(\Gamma^j),$$
and the last sum tends to zero as $i\rightarrow\infty$, since they
are the tails of the convergent sum $\sum_{j\geq
0}\nu_1(\Gamma^j)\leq
\nu_{\alpha-\frac{p_1(d-1)}{q_1(d+1)}}(\Gamma)<\infty$, by
hypothesis. This case follows the previous one and the sufficient
condition on $\alpha$ for (\ref{e:Tmlkv_T-L_Shrlts}) to hold is
proved.

Suppose now that (\ref{e:Tmlkv_T-L_Shrlts}) holds. Fix $j\geq 0$,
$\abs{[\ell]}\preceq2^j$, and $N\in\mathbb{N}$. Let
$\gamma=q_1(s_2-s_1-\frac{1}{p_2})$. Consider the set
$\Gamma_N^{j,\ell}=\{Q_{j,\ell,k}(x):k=(k_1,k_2,\ldots,k_d), 0\leq
k_1,k_2,\ldots,k_d<N \}$ of $N^d$ disjoint anisotropic
parallelograms of area $\abs{Q_j}=2^{-j(d+1)}$. On one hand we have
that
$$\norm{\sum_{P\in\Gamma_N^{j\ell}}\frac{\mathbf{e}_P}{\norm{\mathbf{e}_P}_{\mathbf{f}_2}}}_{\mathbf{f}_1}
    = \left(\int_{\mathbb{R}^d} \left[S^\gamma_{\Gamma_N^{j,\ell}}(x)\right]^{p_1/q_1} dx\right)^{1/p_1}.$$
Since
$$S^\gamma_{\Gamma_N^{j,\ell}}(x)=(2^{-j(d+1)})^\gamma\sum_{P_j\in\Gamma_N^{j,\ell}}\chi_{P_j}(x)=(2^{-j(d+1)})^\gamma\chi_{A^{-j}B^{-\ell}(NP)}(x),$$
then,
$$\norm{\sum_{P\in\Gamma_N^{j\ell}}\frac{\mathbf{e}_P}{\norm{\mathbf{e}_P}_{\mathbf{f}_2}}}_{\mathbf{f}_1}
    = (2^{-j(d+1)})^{\gamma/q_1}[N^d(2^{-j(d+1)})]^{1/p_1}.$$
On the other hand, we have from the right hand side of
(\ref{e:Tmlkv_T-L_Shrlts}) that
$$\nu_{\alpha-\frac{p_1(d-1)}{q_1(d+1)}}(\Gamma_N^{j,\ell})
    =\sum_{P\in\Gamma_N^{j,\ell}}\abs{P}^{\alpha-\frac{p_1(d-1)}{q_1(d+1)}}
    = N^d2^{-j(d+1)(\alpha-\frac{p_1(d-1)}{q_1(d+1)})}.$$
Then, by hypothesis
$$(2^{-j(d+1)})^{\frac{\gamma}{q_1}}(2^{-j(d+1)})^{\frac{1}{p_1}}\lesssim (2^{-j(d+1)})^{\frac{1}{p_1}(\alpha-\frac{p_1(d-1)}{q_1(d+1)})}.$$
From this we deduce $\frac{\gamma}{q_1}+\frac{1}{p_1}\geq
\frac{\alpha}{p_1}-\frac{(d-1)}{q_1(d+1)}$, which implies
$$p_1(s_2-\frac{1}{p_2}-s_1+\frac{1}{p_1}) +
\frac{p_1(d-1)}{q_1(d+1)}\geq \alpha.$$ Similarly, from the left
hand side of (\ref{e:Tmlkv_T-L_Shrlts}) we obtain
$\frac{1}{p_1}(\alpha+\frac{d-1}{d+1})\geq
\frac{\gamma}{q_1}+\frac{1}{p_1},$ which implies $$\alpha\geq
p_1(s_2-\frac{1}{p_2}-s_1+\frac{1}{p_1})-\frac{d-1}{d+1},$$ and the
proof is complete.

\hfill $\blacksquare$ \vskip .5cm   

The second main result is:
\begin{thm}\label{t:Tmlkv_Bsv_Shrlts}
Let $s_1,s_2\in \mathbb{R}$, $0<p_1,p_2<\infty$, $0<q_2\leq\infty$
and
$$\mathbf{u}=\{\norm{\mathbf{e}_P}_{\mathbf{b}^{s_2,q_2}_{p_2}(AB)}\}_{P\in\mathcal{Q}_{AB}}=\{\abs{P}^{-s_2+\frac{1}{p_2}-\frac{1}{2}}\}_{P\in\mathcal{Q}_{AB}}.$$
Then,
\begin{equation}\label{e:Tmlkv_Bsv_Shrlts}
\norm{\sum_{P\in\Gamma}\frac{\mathbf{e}_P}{\norm{\mathbf{e}_P}_{\mathbf{b}^{s_2,q_2}_{p_2}(AB)}}}_{\mathbf{b}^{s_1,p_1}_{p_1}(AB)}
    =\left(\nu_\alpha(\Gamma)\right)^{1/p_1},
\end{equation}
for all $\Gamma\subset \mathcal{Q}_{AB}$ such that
$\nu_{\alpha}(\Gamma)<\infty$ if and only if
$\alpha=p_1[s_2-\frac{1}{p_2}-s_1+\frac{1}{p_1}]$.
\end{thm}
\textbf{Proof}. Consider again the simplified notation. It is a
straight consequence of the definitions
$$\norm{\sum_{P\in\Gamma}\frac{\mathbf{e}_P}{\norm{\mathbf{e}_P}_{\mathbf{b}_2}}}_{\mathbf{b}_1}
    = \left(\sum_{P\in\Gamma}[\abs{P}^{s_2-s_1+\frac{1}{p_1}-\frac{1}{p_2}}]^{p_1}\right)^{1/p_1}
    = \left(\sum_{P\in\Gamma}
    \abs{P}^\alpha\right)^{1/p_1}=\left(\nu_\alpha(\Gamma)\right)^{1/p_1}.
$$
\hfill $\blacksquare$ \vskip .5cm   

This result also applies to shear anisotropic inhomogeneous
Triebel-Lizorkin spaces since
$\mathbf{f}^{s_1,p_1}_{p_1}(AB)=\mathbf{b}^{s_1,p_1}_{p_1}(AB)$.

\vskip0.5cm
\section{Approximation and interpolation}\label{S:RNLApprx_SeqSpcsShrlts}
Once we have found the lower and upper democracy bounds of the
spaces $\mathbf{f}^{s,q}_{p}(AB)$ and $\mathbf{b}^{s,p}_{p}(AB)$ in
(\ref{e:Tmlkv_T-L_Shrlts}) and (\ref{e:Tmlkv_Bsv_Shrlts}), we can
apply the results in Section \ref{sS:RNLASS}.

Proofs for Theorems \ref{t:rhs_demcrcy_T-LShrlts} and
\ref{t:lhs_demcrcy_T-LShrlts} are straight applications of Theorem
\ref{t:J-B_ineq_equiv_for_seq} to the right and left-hand sides of
Theorem \ref{t:Tmlkv_T-L_Shrlts}, respectively. Similarly, Theorem
\ref{t:Charactn_p-eq-q} is consequence of a straight application of
Theorem \ref{t:J-B_ineq_equiv_for_seq} to Theorem
\ref{t:Tmlkv_Bsv_Shrlts}. Part ii) in Theorems
\ref{t:rhs_demcrcy_T-LShrlts}, \ref{t:lhs_demcrcy_T-LShrlts},
\ref{t:Charactn_p-eq-q} and \ref{t:Interpol_LrtzSpcs_Shrlts} is the
special case when Lorenz spaces can be identified with spaces
$\mathbf{b}^{s,p}_{p}(AB)=\mathbf{f}^{s,p}_{p}(AB)$ via Lemma
\ref{l:BsvSeqSpcs_r_LrntzSpcs}.

\vskip0.3cm
\begin{thm}\label{t:rhs_demcrcy_T-LShrlts}
For $s_1,s_2\in\mathbb{R}$, $0<p_1,p_2<\infty$,
$0<q_1,q_2\leq\infty$, $q_1\neq p_1$ with
$\mathbf{u}=\{\norm{\mathbf{e}_P}_{f^{s_2,q_2}_{p_2}(AB)}\}_{P\in\mathcal{Q}_{AB}}$,
fix $\alpha=p_1(s_2-\frac{1}{p_2}-s_1+\frac{1}{p_1})$,
$\beta=\alpha-\frac{p_1(d-1)}{q_1(d+1)}$,
$\xi\in(0,\infty)$ and $\mu\in(0,\infty]$.
\begin{itemize}
  \item[i)] Let $r$ be such that $1/r=\xi+1/p_1$. The following are equivalent:
    \begin{itemize}
      \item[a)] There exists $C>0$ such that for all
        $\Gamma\subset\mathcal{Q}_{AB}$ with
        $\nu_{\beta}(\Gamma)<\infty$,
        $$\norm{\sum_{P\in\Gamma}\frac{\mathbf{e}_P}{u_P}}_{\mathbf{f}^{s_1,q_1}_{p_1}(AB)}
            \leq C \left(\nu_{\beta}(\Gamma)\right)^{1/p_1}.$$
      \item[b)] $\ell^{r,\mu}(\mathbf{u},\nu_{\beta},\mathcal{Q}_{AB})
            \hookrightarrow \mathcal{A}_\mu^\xi(\mathbf{f}^{s_1,q_1}_{p_1}(AB), \nu_{\beta})$.
      \item[c)] The space
            $\ell^{r,\mu}(\mathbf{u},\nu_{\beta},\mathcal{Q}_{AB})$
            satisfies Jackson's inequality of order $\xi$:
            There exists $C>0$ such that
                $$\sigma_{\nu_{\beta}}(t,\mathbf{s})_{\mathbf{f}^{s_1,q_1}_{p_1}(AB)}
                \leq C t^{-\xi}\norm{\mathbf{s}}_{\ell^{r,\mu}(\mathbf{u},\nu_{\beta},\mathcal{Q}_{AB})}.$$
    \end{itemize}

  \item[ii)] If, additionally, $\mu=r$
        and $\gamma=s_1+\frac{d-1}{q_1(d+1)}+\xi[1-\beta]$, the
        following are equivalent:
        \begin{itemize}
            \item[a)] There exists $C>0$ such that for all
                $\Gamma\in\mathcal{Q}_{AB}$ with
                $\nu_{\beta}(\Gamma)<\infty$,
                $$\norm{\sum_{P\in\Gamma}
                    \frac{\mathbf{e}_P}{u_P}}_{\mathbf{f}^{s_1,q_1}_{p_1}(AB)}
                    \leq C \left(\nu_{\beta}(\Gamma)\right)^{1/p_1}.$$
            \item[b)] $\mathbf{b}^{\gamma,r}_r(AB)=\ell^r(\mathbf{u},\nu_{\beta},\mathcal{Q}_{AB})
                \hookrightarrow \mathcal{A}_r^\xi(\mathbf{f}^{s_1,q_1}_{p_1}(AB),
                \nu_{\beta})$.
            \item[c)] The space
                $\mathbf{b}^{\gamma,r}_r(AB)=\ell^r(\mathbf{u},\nu_{\beta},\mathcal{Q}_{AB})$
                satisfies Jackson's inequality of order
                $\xi$:
                There exists $C>0$ such that
                $$\sigma_{\nu_{\beta}}(t,\mathbf{s})_{\mathbf{f}^{s_1,q_1}_{p_1}(AB)}
                \leq
                C t^{-\xi}\norm{\mathbf{s}}_{\mathbf{b}^{\gamma,r}_r(AB)}.$$
        \end{itemize}
\end{itemize}
\end{thm}

\vskip0.3cm
\begin{thm}\label{t:lhs_demcrcy_T-LShrlts}
For $s_1,s_2\in\mathbb{R}$, $0<p_1,p_2<\infty$,
$0<q_1,q_2\leq\infty$, $q_1\neq p_1$ with
$\mathbf{u}=\{\norm{\mathbf{e}_P}_{f^{s_2,q_2}_{p_2}(AB)}\}_{P\in\mathcal{Q}_{AB}}$,
fix $\alpha=p_1(s_2-\frac{1}{p_2}-s_1+\frac{1}{p_1})$,
$\beta=\alpha+\frac{d-1}{d+1}$,
$\xi\in(0,\infty)$ and $\mu\in(0,\infty]$.
\begin{itemize}
  \item[i)] Let $r$ be such that $1/r=\xi+1/p_1$. The following are equivalent:
    \begin{itemize}
        \item[a)] There exists $C>0$ such that for all
            $\Gamma\in\mathcal{Q}_{AB}$ with
            $\nu_{\beta}(\Gamma)<\infty$,
            $$C\left(\nu_{\beta}(\Gamma)\right)^{1/p_1}
                \leq \norm{\sum_{P\in\Gamma}
                \frac{\mathbf{e}_P}{u_P}}_{\mathbf{f}^{s_1,q_1}_{p_1}(AB)}.$$
        \item[b)] $\mathcal{A}^\xi_\mu(\mathbf{f}^{s_1,q_1}_{p_1}(AB),\nu_{\beta})
            \hookrightarrow
            \ell^{r,\mu}(\mathbf{u},\nu_{\beta},\mathcal{Q}_{AB})$.
        \item[c)] The space
            $\ell^{r,\mu}(\mathbf{u},\nu_{\beta},\mathcal{Q}_{AB})$
            satisfies Bernstein's inequality of order $\xi$: There
            exists $C>0$ such that
            $$\norm{\mathbf{s}}_{\ell^{r,\mu}(\mathbf{u},\nu_{\beta},\mathcal{Q}_{AB})}
                \leq C t^\xi
                \norm{\mathbf{s}}_{\mathbf{f}^{s_1,q_1}_{p_1}(AB)},$$
            for all
            $\mathbf{s}\in\Sigma_{t,\nu_{\beta}}\cap\mathbf{f}^{s_1,q_1}_{p_1}(AB)$.
    \end{itemize}

  \item[ii)] If, additionally, $\mu=r$
        and $\gamma=s_1-\frac{d-1}{p_1(d+1)}+\xi(1-\beta)$, the
        following are equivalent:
    \begin{itemize}
        \item[a)] There exists $C>0$ such that for all
            $\Gamma\in\mathcal{Q}_{AB}$ with
            $\nu_{\beta}(\Gamma)<\infty$,
            $$C\left(\nu_{\beta}(\Gamma)\right)^{1/p_1}
                \leq \norm{\sum_{P\in\Gamma}
                \frac{\mathbf{e}_P}{u_I}}_{\mathbf{f}^{s_1,q_1}_{p_1}(AB)}.$$
        \item[b)] $\mathcal{A}^\xi_r(\mathbf{f}^{s_1,q_1}_{p_1}(AB),\nu_{\beta})
            \hookrightarrow
            \ell^r(\mathbf{u},\nu_{\beta},\mathcal{Q}_{AB})=\mathbf{b}^{\gamma,r}_r(AB)$.
        \item[c)] The space
            $\mathbf{b}^{\gamma,r}_r(AB)=\ell^r(\mathbf{u},\nu_{\beta},\mathcal{Q}_{AB})$
            satisfies Bernstein's inequality of order $\xi$: There
            exists $C>0$ such that
            $$\norm{\mathbf{s}}_{\mathbf{b}^{\gamma,r}_r(AB)}
                \leq C t^\xi
                \norm{\mathbf{s}}_{\mathbf{f}^{s_1,q_1}_{p_1}(AB)},$$
            for all $\mathbf{s}\in\Sigma_{t,\nu_{\beta}}\cap\mathbf{f}^{s_1,q_1}_{p_1}(AB).$
    \end{itemize}
\end{itemize}
\end{thm}

\vskip0.3cm

In the case $q_1=p_1<\infty$, \emph{i.e.}
$\mathbf{b}^{s_1,p_1}_{p_1}(AB)=\mathbf{f}^{s_1,p_1}_{p_1}(AB)$,
applying Theorem \ref{t:J-B_ineq_equiv_for_seq} to Theorem
\ref{t:Tmlkv_Bsv_Shrlts} yields the next characterizations.
\begin{thm}\label{t:Charactn_p-eq-q}
For $s_1,s_2\in\mathbb{R}$, $0<p_1,p_2\leq\infty$ with
$\mathbf{u}=\{\norm{\mathbf{e}_Q}_{b^{s_2,q_2}_{p_2}(AB)}\}_{Q\in\mathcal{Q}_{AB}}$,
fix $\alpha=p_1(s_2-\frac{1}{p_2}-s_1+\frac{1}{p_1})$,
$\xi\in(0,\infty)$ and $\mu\in(0,\infty]$.
\begin{itemize}
    \item[i)] Let $r$ be such that $1/r=\xi+1/p_1$. The following are equivalent:
    \begin{itemize}
        \item[a)] The Temlyakov property holds:
            $$\norm{\sum_{Q\in\Gamma}\frac{\mathbf{e}_Q}{u_Q}}_{\mathbf{b}^{s_1,p_1}_{p_1}(AB)}=\left(\nu_\alpha(\Gamma)\right)^{1/p_1}.$$
        \item[b)]
            $\ell^{r,\mu}(\mathbf{u}, \nu_\alpha,\mathcal{Q}_{AB})=\mathcal{A}_\mu^\xi(\mathbf{b}^{s_1,p_1}_{p_1}(AB),\nu_\alpha).                  $
        \item[c)] The Jackson's inequality of order $\xi$ holds: There exists
            $C>0$ such that
            $$\sigma_{\nu_\alpha}(t,\mathbf{s})_{\mathbf{b}^{s_1,p_1}_{p_1}(AB)}
                \leq Ct^{-\xi}\norm{\mathbf{s}}_{\ell^{r,\mu}(\mathbf{u}, \nu_\alpha,\mathcal{Q}_{AB})}.$$
            The Bernstein's inequality of order $\xi$ holds: There exists
            $C'>0$ such that
            $$\norm{\mathbf{s}}_{\ell^{r,\mu}(\mathbf{u}, \nu_\alpha,\mathcal{Q}_{AB})}
                \leq
                C't^\xi\norm{\mathbf{s}}_{\mathbf{b}^{s_1,p_1}_{p_1}(AB)},$$
            for all $\mathbf{s}\in\Sigma_{t,\nu_{\alpha}}\cap\mathbf{b}^{s_1,p_1}_{p_1}(AB).$
    \end{itemize}
    \item[ii)] If, additionally, $\mu=r$ and $\gamma=s_1+\xi(1-\alpha)$, the following are
        equivalent:
    \begin{itemize}
        \item[a)] The Temlyakov property holds:
            $$\norm{\sum_{Q\in\Gamma}\frac{\mathbf{e}_Q}{u_Q}}_{\mathbf{b}^{s_1,p_1}_{p_1}(AB)}=\left(\nu_\alpha(\Gamma)\right)^{1/p_1}.$$
        \item[b)] $\mathbf{b}^{\gamma,r}_r(AB)=\ell^{r}(\mathbf{u}, \nu_\alpha,\mathcal{Q}_{AB})
                =\mathcal{A}_r^\xi(\mathbf{b}^{s_1,p_1}_{p_1}(AB),\nu_\alpha).$
        \item[c)] The Jackson's inequality of order $\xi$ holds: There exists
            $C>0$ such that
            $$\sigma_{\nu_\alpha}(t,\mathbf{s})_{\mathbf{b}^{s_1,p_1}_{p_1}(AB)}
                \leq Ct^{-\xi}\norm{\mathbf{s}}_{\mathbf{b}^{\gamma,r}_r(AB)}.$$
            The Bernstein's inequality of order $\xi$ holds: There exists
            $C'>0$ such that
            $$\norm{\mathbf{s}}_{\mathbf{b}^{\gamma,r}_r(AB)}
                \leq
                C't^\xi\norm{\mathbf{s}}_{\mathbf{b}^{s_1,p_1}_{p_1}(AB)},$$
            for all $\mathbf{s}\in\Sigma_{t,\nu_{\alpha}}\cap\mathbf{b}^{s_1,p_1}_{p_1}(AB).$
    \end{itemize}
    \end{itemize}
\end{thm}

\vskip0.3mm
When $0<p_1,p_2<\infty$, Theorem (\ref{t:Charactn_p-eq-q}) also applies to
$\mathbf{f}^{s_i,p_i}_{p_i}, i=1,2$.

We finish with a result on interpolation of shear anisotropic
inhomogeneous Besov spaces.

\begin{thm}\label{t:Interpol_LrtzSpcs_Shrlts}
For $s_1,s_2\in\mathbb{R}$, $0<p_1,p_2\leq\infty$ and
        $\alpha=p_1(s_2-\frac{1}{p_2}-s_1+\frac{1}{p_1})$ with
        $\mathbf{u}=\{\abs{Q}^{-s_2+\frac{1}{p_2}-\frac{1}{2}}\}_{Q\in\mathcal{Q}_{AB}}$,
        fix $\xi_i\in(0,\infty)$ and $\mu_i\in(0,\infty]$ and let $r_i$
        be such that $1/r_i=\xi_i+1/p_1$ for $i=0,1$. Then,
\begin{enumerate}
  \item[i)] For all
        $\theta\in(0,1)$ and $\mu\in(0,\infty]$,
        $$(\mathcal{A}_{\mu_0}^{\xi_0}(\mathbf{b}^{s_1,p_1}_{p_1}(AB),\nu_\alpha),
        \mathcal{A}_{\mu_1}^{\xi_1}(\mathbf{b}^{s_1,p_1}_{p_1}(AB),\nu_\alpha))_{\theta,\mu}
        =\mathcal{A}_\mu^\xi(\mathbf{b}^{s_1,p_1}_{p_1}(AB),\nu_\alpha),$$
        where $\xi=(1-\theta)\xi_0+\theta\xi_1$ and, therefore,
        $$(\ell^{r_0,\mu_0}(\mathbf{u},\nu_\alpha,\mathcal{Q}_{AB}),
        \ell^{r_1,\mu_1}(\mathbf{u},\nu_\alpha,\mathcal{Q}_{AB}))_{\theta,\mu}
        =\ell^{r,\mu}(\mathbf{u},\nu_\alpha,\mathcal{Q}_{AB}),$$
        where $1/r=(1-\theta)/r_0+\theta/r_1=(1-\theta)(\xi_0+\frac{1}{p_1})+\theta(\xi_1+\frac{1}{p_1})=\xi+\frac{1}{p_1}$.
  \item[ii)] If, additionally, $\mu_i=r_i$ and $\mu=r$, let $\gamma_0=s_1+\xi_0(1-\alpha)$ and
        $\gamma_1=s_1+\xi_1(1-\alpha)$. Then, for all $\theta\in(0,1)$
        and $r\in(0,\infty)$,
        $$(\mathbf{b}^{\gamma_0,r_0}_{r_0}(AB),
            \mathbf{b}^{\gamma_1,r_1}_{r_1}(AB))_{\theta,r}
            =\mathbf{b}^{\gamma,r}_r(AB),$$
        where
        $\gamma=(1-\theta)\gamma_0+\theta\gamma_1=s_1+(1-\theta)(\xi_0(1-\alpha))+\theta(\xi_1(1-\alpha))
            =s_1+\xi(1-\alpha)$.
\end{enumerate}
\end{thm}
When $0<p_1,p_2<\infty$, Theorem (\ref{t:Interpol_LrtzSpcs_Shrlts}) also applies to
$f^{s_i,p_i}_{p_i}, i=1,2$.
\begin{rem}\label{r:t:Interpol_LrtzSpcs_Shrlts}
It is not hard to see that one can choose $s_1,s_2,p_1,p_2$ in
Theorem \ref{t:Interpol_LrtzSpcs_Shrlts} such that we have
interpolation of shear anisotropic inhomogeneous Besov spaces in
point ii) of Theorem \ref{t:Interpol_LrtzSpcs_Shrlts} for any
$r_i\in(0,\infty)$, $\gamma_i=s_1+\xi_i(1-\alpha)\in\mathbb{R}$,
$\theta\in(0,1)$ and $r\in(0,\infty)$, since $r_i=r_i(\xi_i,p_1)$,
$\mathbf{u}=\mathbf{u}(s_2,p_2)$ and
$\alpha=\alpha(s_1,p_1,s_2,p_2)$ can be chosen independently.
Observe also that $\gamma$ in point ii) of Theorem
\ref{t:Interpol_LrtzSpcs_Shrlts} coincides with that in point ii) of
Theorem \ref{t:Charactn_p-eq-q}.
\end{rem}

\section{Comments}\label{S:Extn_Parblc_Molecls}

\subsection{Extension to parabolic molecules}\label{sS:Extn_ParMol}
Here we extend the comments made in the Introduction, Section
\ref{sS:Intro_Shrlts-SmthnssSpcs}. With Theorem 2.11 in
\cite{BoNi07}, on equivalent admissible coverings implies equivalent
decomposition spaces, it was proved in \cite{LMN2012} that the
shearlet smoothness spaces $S^\beta_{p,q}$ and the curvelet (first
and second generation) decomposition spaces $G^\beta_{p,q}$
\cite{BoNi07} are equivalent with equivalent norms. Since the shear
anisotropic inhomogeneous Besov sequence spaces
$\mathbf{b}^{\alpha,q}_p$ are based on the same \textbf{bounded
admissible partition of unity} (\textbf{BAPU}) as $S^\beta_{p,q}$
(around the cartesian coronae concentrated in
$\abs{\xi}\sim2^{2j}$), Theorems \ref{t:Charactn_p-eq-q} and
\ref{t:Interpol_LrtzSpcs_Shrlts} are immediately extended to
$S^\beta_{p,p}$ and $G^\beta_{p,p}$, previous normalization in the
smoothness parameters that we explain next. When we choose a weight
$2^{j(d+1)\alpha}$ instead of $2^{j\beta}$ (as in \cite{BoNi07} and
in \cite{LMN2012}) the space $S^\beta_{p,q}$ (and therefore $G^\beta_{p,q}$)
coincide with the space $\mathbf{b}^{\alpha,q}_p(AB)$ whose natural
``weights" are $\abs{Q}^{-\alpha}=2^{j(d+1)\alpha}$, \emph{i.e.} the
volumes of the ``cubes" in (\ref{e:def_Qjlk}). This result is only
for band-limited shearlet generators. More recently, Grohs and
Kutyniok \cite{GrKu12} proved the equivalence of more general
parabolic molecules generated sequence spaces and their
approximation properties whenever the parabolic molecules are
sufficiently smooth as well as localized in space and frequency.
Hence, all of our results extend to all spaces generated from
parabolic molecules that are equivalent to the shear anisotropic
inhomogeneous Besov and Tribel-Lizorkin spaces as defined in
\cite{Ver13} and \cite{Ver12} whenever they are sufficiently 
smooth as well as localized in space and frequency.

\subsection{Shearlet spaces and $BV(\mathbb{R}^2)$}\label{sS:BV_Shrlts}
As mentioned in the Introduction, the inclusions
$b^{1,1}_1(\mathcal{D}_+)\hookrightarrow
bv(\mathcal{D}_+)\hookrightarrow w\ell^1(\mathcal{D}_+)$
were proved in \cite{CDDD} for the non-homogeneous basis indexed by $\mathcal{D}_+$.
We can apply the method of retracts to transfer results
in sequence spaces to function/distribution spaces weather we are dealing with wavelets or shearlets.  
By identifying a classical Besov sequence space $b^{s,1}_1(\mathcal{D}_+)$, $s=1-\varepsilon<1$, 
with a weighted discrete Lorenz space $\ell^1(\abs{Q}^{\varepsilon/2})$ (as defined in
Section \ref{sS:RNLASS} with weight sequence $\mathbf{u}=\{\abs{Q}^{\varepsilon/2}\}_{Q\in\mathcal{D}}$)
via Lemma 2.10.4 in \cite{HV11} and from the fact that 
$w\ell^1\hookrightarrow \ell^1(\abs{Q}^{\varepsilon/2})$ we have
$b^{1,1}_1(\mathcal{D}_+)\hookrightarrow
bv(\mathcal{D}_+)\hookrightarrow b^{s,1}_1(\mathcal{D}_+), s=1-\varepsilon<1$.
From the results in \cite{CDH} or
\cite{HV11} one has for $N$-term approximation $\mathcal{A}^\xi_1(L^2)=B^{\gamma,1}_1$.
Therefore, one conjectures that the approximation error measured in
$L^2$ of a $f\in BV$ decays (apply Theorem 2.1 with 
$B^{\gamma,1}_1, \gamma=1,1-\varepsilon$ to obtain Jackson's inequality) as $N^{-\xi}$
with $\xi\in[(1-\varepsilon)/2,1/2]$. Nevertheless, the sharp value $\xi=1/2$ is not
proved by the results in neither \cite{CDH} nor \cite{HV11} but by
the more thorough arguments in \cite{CDDD}. One is tempted to think
that shearlets might be equally as good as wavelets to approximate a
function in $BV(\mathbb{R}^2)$. However, the shearlet representation
is more redundant than that of wavelets since it involves
directionality. Let us show what we can get from our results. From
Theorem 4.2 in \cite{Ver13} one has
$B^{1-\varepsilon,1}_1(\mathbb{R}^2)\hookrightarrow
\mathbf{B}^{\gamma,1}_1(AB)$ for
$\gamma<-\frac{1}{3}(1+2\varepsilon)$, and from Theorem 4.3 in
\cite{Ver13} one has $\mathbf{B}^{s,2}_2(AB)\hookrightarrow
L^2(\mathbb{R}^2)$ for $s>1$. 
Hence, ii) in Theorem \ref{t:Charactn_p-eq-q} 
shows that the inequality
$$\sigma_{\nu_\alpha}(t,\mathbf{s})_{\mathbf{b}^{s,2}_2(AB)}\lesssim t^{-\xi}\norm{\mathbf{s}}_{\mathbf{b}^{\gamma,1}_1(AB)},$$
with $\xi=\frac{1}{2}$ and $\gamma$ as above, can only occur if we
do not impose at the same time $N$-term approximation ($\alpha=0$, and hence $t=N$)
and a norm of the approximation error comparable to the $L^2$ norm
($s>1$).



\end{document}